\documentclass[12pt]{article}

\usepackage{a4,graphicx,amsmath,amsfonts,amssymb,bbm,amsthm}
\usepackage{color,url}
\usepackage{mathtools}

\newcommand{\R}{\mathbbm{R}}
\newcommand{\C}{\mathbbm{C}}
\newcommand{\N}{\mathbbm{N}}

\newcommand{\ltwo}{\mathcal{L}^2(\mathcal{M},\rho)}
\newcommand{\ltwotime}{\mathcal{L}^2[0,\infty)}
\newcommand{\htwo}{\mathcal{H}_2}
\newcommand{\hinf}{\mathcal{H}_{\infty}}

\newtheorem{definition}{Definition}
\newtheorem{theorem}{Theorem}

\newtheorem{lemma}{Lemma}

\begin{document}
\begin{center}
  {\bf \Large Stability-preserving model order reduction \\[1ex]
    for linear stochastic Galerkin systems}

  \vspace{10mm}

{\large Roland~Pulch}

\vspace{0.1cm}

\begin{small}

{Institute of Mathematics and Computer Science, 
University of Greifswald, \\
Walther-Rathenau-Str.~47, 
D-17489 Greifswald, Germany.}\\
Email: {\tt roland.pulch@uni-greifswald.de}

\end{small}
  
\end{center}

\bigskip\bigskip


\begin{center}
{Abstract}

\begin{tabular}{p{13cm}}
  Mathematical modeling often yields linear dynamical systems in
  science and engineering.
  We change physical parameters of the system into random variables 
  to perform an uncertainty quantification.
  The stochastic Galerkin method yields a larger linear dynamical system,
  whose solution represents an approximation of random processes.
  A model order reduction (MOR) of the Galerkin system is advantageous
  due to the high dimensionality.
  However, asymptotic stability may be lost in some MOR techniques.
  In Galerkin-type MOR methods, the stability can be guaranteed by a
  transformation to a dissipative form.
  Either the original dynamical system or the stochastic Galerkin
  system can be transformed.
  We investigate the two variants of this stability-preserving
  approach.
  Both techniques are feasible, while featuring different properties
  in numerical methods. 
  Results of numerical computations are demonstrated for two test examples
  modeling a mechanical application and an electric circuit, respectively.
  
\bigskip


\bigskip

\end{tabular}
\end{center}


\thispagestyle{plain}

\clearpage

\markright{R.~Pulch: Stability-preserving MOR for linear stochastic Galerkin systems}


\section{Introduction}
Numerical simulation of mathematical models represents the main issue
in scientific computing.
We consider linear dynamical systems, which play an important role in
mechanics and electrical engineering, for example.
Furthermore, uncertainty quantification becomes more and more relevant
in many fields of applications, see~\cite{maten-nanocops}, for instance.
A common approach is to replace uncertain parameters by random variables,
see~\cite{sullivan,xiu-book}.
Statistics of the stochastic model can be computed by sampling methods
or quadrature rules.
Alternatively, the stochastic Galerkin method changes the
random-dependent linear dynamical system into a larger deterministic
linear dynamical system. 

The dimension of the stochastic Galerkin system becomes huge in the
case of large numbers of random variables.
Methods of model order reduction (MOR) are able to decrease the
complexity.
Transient solutions of a reduced system allow for an efficient
numerical simulation.
Several MOR methods are available for general linear dynamical systems, 
see \cite{antoulas,benner,benner-mehrmann,schilders}.
MOR of linear stochastic Galerkin systems was also examined
in several previous works
\cite{freitas,mi-etal,pulch-maten,pulch-scee2014,pulch17,zou-etal}.
MOR of nonlinear stochastic Galerkin systems was considered
in~\cite{pulch19}.

However, even though the linear stochastic Galerkin system is
asymptotically stable, the reduced Galerkin system often looses
this stability in some MOR techniques.
We investigate stability-preserving strategies in the case of
Galerkin-type projection-based MOR like the Arnoldi method or
proper orthogonal decomposition, for example.
Galerkin-type MOR can be applied to any linear dynamical system
(not only stochastic Galerkin systems). 
A dissipativity property guarantees the preservation of stability.
If a general linear dynamical system does not satisfy
the dissipativity property,
then it can be transformed into a dissipative structure, 
see~\cite{castane-selga,prajna,pulch-naco}.
The crucial part to identify a transformation consists in the solution
of a Lyapunov equation.
Direct methods, see~\cite{hammarling}, or approximate methods,
see~\cite{penzl1998,penzl-sisc,son-stykel,wolf},
yield the numerical solutions of Lyapunov equations.

We examine the stability-preserving approach in the case of linear stochastic
Galerkin systems consisting of ordinary differential equations.
Several variants are feasible.
The high-dimensional Galerkin-projected system is transformed or,
vice versa, the original systems are transformed followed by a
Galerkin projection.
We analyze the two strategies and another variant.

In addition, network approaches produce models consisting of 
differential-alge\-braic equations in industrial applications.
Thus we extend the stability-preserving techniques to this class
of problems.
The Lyapunov equations have no solution now.
Therefore, we use a regularization technique, which was also employed
in~\cite{mueller}.

We apply the analyzed techniques to mathematical models of two test examples:
a mass-spring-damper system and an electric circuit of a band-pass filter.


\section{Stability preservation in reduction}
\label{sec:mor}
We review a concept for stability preservation in Galerkin-type
projection-based MOR for general linear dynamical systems.

\subsection{Linear dynamical systems}
We consider linear dynamical systems in the form
\begin{align} \label{ode}
  \begin{split}
    E \dot{x}(t) & = A x(t) + B u(t) \\[1ex]
    y(t) & = C x(t) 
  \end{split}
\end{align}
with constant matrices $A,E \in \R^{n \times n}$, $B \in \R^{n \times n_{\rm in}}$,
and $C \in \R^{n_{\rm out} \times n}$.
The state variables or inner variables are
$x: [0,t_{\rm end}] \rightarrow \R^n$.
Inputs $u: [0,t_{\rm end}] \rightarrow \R^{n_{\rm in}}$
are supplied to the system.
The outputs 
$y: [0,t_{\rm end}] \rightarrow \R^{n_{\rm out}}$
are defined as quantities of interest (QoI).
Initial value problems are given by $x(0) = x_0$.

If the mass matrix~$E$ is non-singular, then
the system~(\ref{ode}) consists of ordinary differential equations (ODEs).
If the mass matrix~$E$ is singular,
then the system~(\ref{ode}) represents
differential-algebraic equations (DAEs).
Furthermore, we assume that the system satisfies the following
stability condition.

We specify some common notions.

\begin{definition} \label{def:abscissa}
  Given a matrix pencil $(E,A)$ with $E,A \in \R^{n \times n}$,
  the set of {\em eigenvalues} is 
  $\Lambda = \left\{ \lambda \in \C \, : \, \det( \lambda E - A) = 0 \right\}$
  and the {\em spectral abscissa} reads as
  $\alpha(E,A) = \max
  \left\{ {\rm Re}(\lambda) \, : \, \lambda \in \Lambda \right\}$.
  The {\em spectral abscissa} of a single matrix~$A$ is the
  spectral abscissa of the matrix pencil $(I,A)$ with the identity matrix
  $I \in \R^{n \times n}$.
\end{definition}

\begin{definition} \label{def:matrix-pencil}
  A matrix pencil $(E,A)$ with $E,A \in \R^{n \times n}$ is called 
  {\em regular}, if there is (at least one) $\lambda \in \C$ such that
  $\det( \lambda E - A) \neq 0$. 
\end{definition}

\begin{definition} \label{def:stable}
  A linear dynamical system~(\ref{ode}) is called
  {\em asymptotically stable}, if the involved matrix pencil $(E,A)$
  has a spectral abscissa satisfying $\alpha(E,A) < 0$.
\end{definition}

The asymptotic stability implies that the generalized eigenvalue problem
has only a finite number of eigenvalues,
which all exhibit a negative real part.
We assume that the system~(\ref{ode}) satisfies this
stability condition.
In the case of ODEs, the matrix pencil is always regular,
even if the system is unstable.
In the case of DAEs, the asymptotic stability implies that the
matrix pencil $(E,A)$ is regular.

\subsection{Transfer functions and Hardy norms}
The input-output behavior of a linear dynamical system~(\ref{ode})
can be described in the frequency domain,
see~\cite{antoulas} for ODEs and~\cite{benner-stykel} for DAEs.

\begin{definition} \label{def:transfer-function}
  The {\em transfer function} of a linear dynamical system~(\ref{ode})
  is the mapping
  $H : \C \backslash \Lambda \rightarrow \C^{n_{\rm out} \times n_{\rm in}}$ with
  \begin{equation} \label{transfer-function}
    H (s) = C (sE-A)^{-1} B ,
  \end{equation}
  where $\Lambda$ represents the set of eigenvalues
  from Definition~\ref{def:abscissa}.
\end{definition}

The transfer function is always a rational function,
whose poles are the eigenvalues.
The regularity of the matrix pencil guarantees a finite set of poles.
The magnitude of a transfer function can be measured by
Hardy norms, see~\cite{antoulas}.

\begin{definition} \label{def:hardy-norm}
  The $\htwo$-{\em norm} of a transfer function reads as
  \begin{equation} \label{htwo-norm}
    \| H \|_{\htwo} = \sqrt{ \frac{1}{2\pi}
      \int_{-\infty}^{+\infty} \| H({\rm i}\omega) \|_{\rm F}^2 \; {\rm d}\omega }
  \end{equation}
  including the Frobenius matrix norm $\| \cdot \|_{\rm F}$,
  the angular frequency~$\omega$, and ${\rm i} = \sqrt{-1}$.
\end{definition}

Since we assume asymptotically stable linear dynamical systems~(\ref{ode}),
the transfer function is defined on the complete imaginary axis. 
In the case of ODEs, the $\htwo$-norm is always finite.
In the case of DAEs, the existence of the $\htwo$-norm is not guaranteed.
Nevertheless, the $\htwo$-norm is quite often finite for DAEs
of index~1 or~2. 
The norm~(\ref{htwo-norm}) may also exist for unstable systems, if there
is no pole on the imaginary axis.


\subsection{Projection-based model order reduction}
Projection matrices $V,W \in \R^{n \times r}$ of full rank are specified
with $r \ll n$.
Concerning the full-order model (FOM) in~(\ref{ode}),
the reduced-order model (ROM) reads as
\begin{align} \label{rom}
  \begin{split} 
    \bar{E} \dot{\bar{x}}(t) & = \bar{A} \bar{x}(t) + \bar{B} u(t) \\[1ex]
    \bar{y}(t) & = \bar{C} \bar{x}(t) 
  \end{split}
\end{align}
with state variables or inner variables
$\bar{x} : [0,t_{\rm end}] \rightarrow \R^r$.
Initial values $\bar{x}(0) = \bar{x}_0$ are supposed.
We obtain the matrices via
\begin{equation} \label{matrices-reduced}
  \bar{A} = W^\top A V , \quad
  \bar{B} = W^\top B , \quad
  \bar{C} = C V , \quad
  \bar{E} = W^\top E V ,
\end{equation}
which is also called a Petrov-Galerkin-type MOR.
A Galerkin-type projection-based MOR is characterized by $W=V$,
where just one projection matrix has to be determined.
Important examples are the one-sided Arnoldi method
and the proper orthogonal decomposition (POD), see~\cite{antoulas}.

Each linear dynamical system is described by a
transfer function~(\ref{transfer-function})
in the frequency domain.
The difference between the transfer functions~$H$ of FOM and
$\bar{H}$ of ROM quantifies the error of the MOR.
Hardy norms like the $\htwo$-norm from Definition~\ref{def:hardy-norm},
for example, can be applied to their transfer functions.
However, a small error in the frequency domain implies a small error
in the time domain only if both systems are asymptotically stable.
It holds that
\begin{equation} \label{error-time-domain}
  \sup_{t \ge 0} \| y(t) - \bar{y}(t) \|_{\infty} \le
  \left\| H - \bar{H} \right\|_{\htwo} \| u \|_{\ltwotime}
\end{equation}
with the maximum vector norm $\| \cdot \|_{\infty}$ and the
$\ltwotime$-norm in time provided that all initial values are zero,
see~\cite{benner-gugercin}.
The $\htwo$-norm of Definition~\ref{def:hardy-norm} is a
strong measure.
Hence there are neither a priori error bounds nor cheap a posteriori
error bounds available in MOR techniques.
Just an approximation of the $\htwo$-norm in~(\ref{error-time-domain})
can be computed posterior, where the computational effort is
dominated by evaluations of the transfer function in the FOM.
The balanced truncation method yields an a priori error bound
in the $\hinf$-norm, see~\cite[p.~212]{antoulas}.

\subsection{Dissipative systems}
In balanced truncation, see~\cite{antoulas}, the ROM~(\ref{rom}) is always
asymptotically stable provided that the system~(\ref{ode}) is
asymptotically stable.
Yet the asymptotic stability may be lost in the ROM~(\ref{rom})
within other MOR methods like Krylov subspace techniques,
see~\cite{freund}, and POD, for example.

Using Galerkin-type MOR, the stability is guaranteed for some classes
of linear dynamical systems.
In the case of ODEs, we define the following type of system.

\begin{definition} \label{def:dissipative}
  A linear dynamical system~(\ref{ode}) is called {\em dissipative}, if
  \begin{enumerate}
  \item $E$ is symmetric as well as positive definite, and
  \item $A+A^\top$ is negative definite.
  \end{enumerate}
\end{definition}

The above condition represents a dissipativity of the matrix~$A$
as shown in~\cite{panzeretal}.
Other definitions of dissipative systems are used in the literature.
We prove a property, which was shown for an explicit system of ODEs
in~\cite{prajna}.

\begin{theorem} \label{thm:dissipative-is-stable}
  If the linear dynamical system~(\ref{ode}) is dissipative
  with respect to Definition~\ref{def:dissipative},
  then it is also asymptotically stable as in Definition~\ref{def:stable}.
\end{theorem}

\underline{Proof:}

Let $E=LL^\top$ be the Cholesky decomposition of the mass matrix.
The system~(\ref{ode}) is equivalent to the explicit system of ODEs
\begin{equation} \label{ode-trafo-thm}
  \dot{z}(t) = L^{-1} A L^{-\top} z(t) + L^{-1} B u(t)
\end{equation}
with $z = L^\top x$ and $y = C L^{-\top} z$.
The symmetric part of the involved system matrix~$\tilde{A}$ reads as
$$  S := \tilde{A} + \tilde{A}^\top =
L^{-1} A L^{-\top} + \left( L^{-1} A L^{-\top} \right)^\top =
L^{-1} \left( A + A^\top \right) L^{-\top} . $$
It holds that $v^* S v = (L^{-\top}v)^* (A+A^\top) (L^{-\top}v)$
for any $v \in \C^n \backslash \{ 0 \}$.
Thus $S$ is negative definite.
Let $\lambda \in \C$ be an eigenvalue of $\tilde{A}$ and $v$ be
an associated eigenvector satisfying $v^* v = 1$.
It follows that
$$ 2 \, {\rm Re}(\lambda) = \lambda + \overline{\lambda} =
\left( \lambda + \overline{\lambda} \right) v^*v =
v^* \left( \tilde{A} + \tilde{A}^\top \right) v = v^* S v < 0 . $$
Thus the systems~(\ref{ode-trafo-thm}) and~(\ref{ode})
are asymptotically stable.
\hfill $\Box$

\medskip

In contrast, the asymptotic stability does not imply the
dissipativity of Definition~\ref{def:dissipative},
even if the mass matrix~$E$ is symmetric and positive definite.
An example will be presented in Section~\ref{sec:msd}.
Now we consider Galerkin-type MOR.

\begin{theorem} \label{thm:dissipative}
  If the linear dynamical system~(\ref{ode}) is dissipative,
  then a Galerkin-type MOR yields a dissipative
  reduced system~(\ref{rom}).
  Hence the reduced system is asymptotically stable.
\end{theorem}

The proof can be found in~\cite{pulch-naco}, for example.

\subsection{Transformations}
\label{sec:transformation}
The asymptotic stability of a linear dynamical system is invariant
with respect to basis transformations.
In contrast, if a system of ODEs is not dissipative, then it can be
converted to an equivalent dissipative form by a basis transformation
in the state space, see~\cite{prajna}. 
Alternatively, a basis transformation is feasible in the image space only,
see~\cite{castane-selga}.
We require a symmetric positive definite solution~$M \in \R^{n \times n}$
of the Lyapunov inequality
\begin{equation} \label{lyapunov-inequality}
  A^\top M E + E^\top M A < 0 ,
\end{equation}
which means that the matrix on the left-hand side
of~(\ref{lyapunov-inequality}) is negative definite.
If the mass matrix~$E$ is non-singular and the
matrix pencil $(E,A)$ satisfies the Definition~\ref{def:stable}
of asymptotic stability,
then an infinite set of solutions~$M$ exists.

We change the system of ODEs~(\ref{ode}) into the equivalent system
\begin{align} \label{ode-diss}
  \begin{split} 
    E^\top M E \dot{x}(t) & = E^\top M A x (t) + E^\top M B u(t) \\[1ex]
    y(t) & = C x(t) . 
  \end{split}
\end{align}
The transformed system exhibits the desired dissipativity property. 

\begin{theorem} \label{thm:transformation}
  If the asymptotically stable linear dynamical system~(\ref{ode})
  has a non-sin\-gu\-lar mass matrix~$E$,
  then the transformed system~(\ref{ode-diss}) with~$M$
  satisfying~(\ref{lyapunov-inequality}) is dissipative
  in view of Definition~\ref{def:dissipative}.
\end{theorem}

We use a Galerkin-type MOR with a projection matrix~$V$
to the transformed system~(\ref{ode-diss}).
This approach can be written as a Petrov-Galerkin-type MOR
applied to the original system~(\ref{ode})
with matrices~(\ref{matrices-reduced}) and the projection matrix
\begin{equation} \label{projection-w}
  W = M E V .
\end{equation}
Thus we do not require to calculate the transformed system~(\ref{ode-diss})
explicitly.
Alternatively, we compute the projection matrix~(\ref{projection-w}).
However, the original Galerkin-type MOR of the system~(\ref{ode}) is
not equivalent to the MOR of the system~(\ref{ode-diss}).

\subsection{Numerical solution of Lyapunov inequality}
\label{sec:solve-lyapunov}
We solve the Lyapunov inequality~(\ref{lyapunov-inequality}) using
a Lyapunov equation
\begin{equation} \label{lyapunov}
  A^\top M E + E^\top M A + F = 0
\end{equation}
including a predetermined symmetric positive definite matrix
$F \in \R^{n \times n}$.
This matrix represents a degree of freedom, because any choice yields
a symmetric positive definite solution of the Lyapunov inequality.
A simple admissible choice is the identity matrix $I_n \in \R^{n \times n}$.
Moreover, we do not need to solve the Lyapunov
equation~(\ref{lyapunov}) with a high accuracy,
because a rough approximation~$\widetilde{M}$ often still satisfies
the Lyapunov inequality~(\ref{lyapunov-inequality}).
In~\cite{pulch-arxiv}, it was shown that any approximation~$\widetilde{M}$
of the exact solution~$M$ of the Lyapunov equation~(\ref{lyapunov})
for $F=I_n$ with the property
\begin{equation} \label{suff-cond}
  \| \widetilde{M} - M \| <
  \frac{1}{\| A^\top \| \cdot \| E \| + \| A \| \cdot \| E^\top \|}
\end{equation}
in some subordinate matrix norm also solves the
Lyapunov inequality~(\ref{lyapunov-inequality}).
The condition~(\ref{suff-cond}) is just sufficient and not necessary.
However, approximations~$\widetilde{M}$ satisfying~(\ref{suff-cond})
may have a relative error of up to 100\%, see~\cite{pulch-arxiv},
which motivates that rough estimates can solve the problem.

There are direct methods to compute a solution~$M$ of~(\ref{lyapunov})
or a symmetric decomposition $M=L L^\top$,
see~\cite{hammarling,penzl1998}.
Their computational effort is typically $\mathcal{O}(n^3)$.
In the high-dimensional case, we have to use approximate methods
to decrease the computation work.
The following techniques are available:
\begin{itemize}
\item[i)] projection methods (Krylov subspace techniques, POD, etc.),
  see~\cite{kramer-singler,wolf},
\item[ii)] alternating direction implicit (ADI) iteration,
  see~\cite{lu-wachspress,penzl-sisc},
\item[iii)] frequency domain integrals,
  see~\cite{benner-schneider,pulch-arxiv},
\end{itemize}
and others.
In the cases~(i) and~(ii),
the methods yield an approximation $\widetilde{M} = Z Z^\top$
with a low-rank factor $Z \in \R^{n \times k}$ ($k \ll n$).
Thus the transformation is given by a singular matrix~$\widetilde{M}$.
It follows that the mass matrix~$\bar{E}$ of the reduced system~(\ref{rom})
may become singular or ill-conditioned, as shown in~\cite{pulch-naco}.
In contrast, the method~(iii) from~\cite{pulch-arxiv} computes the
projection matrix~(\ref{projection-w}),
where the underlying approximation~$\widetilde{M}$ is always non-singular.
However, the matrix~$\widetilde{M}$ is never computed
but a matrix-matrix product with this approximation.
In the frequency domain integral approach, the projection matrix~$V$
has to be determined by the original linear dynamical
system~(\ref{ode}).


\section{Stochastic Galerkin systems}
\label{sec:galerkin}
We illustrate the concept of the stochastic Galerkin method
and define our problem under investigation.

\subsection{Random linear dynamical systems}
We include parameters in the linear dynamical systems and obtain
\begin{align} \label{ode-parameter}
  \begin{split}
    E(\mu) \dot{x}(t,\mu) & = A(\mu) x(t,\mu) + B(\mu) u(t) \\[1ex]
    y(t,\mu) & = C(\mu) x(t,\mu) .
  \end{split}
\end{align}
The matrices $A,E \in \R^{n \times n}$, $B \in \R^{n \times n_{\rm in}}$, and
$C \in \R^{n_{\rm out} \times n}$
depend on parameters $\mu \in \mathcal{M} \subseteq \R^q$.
We assume that the dimension~$n$ is independent of the number of
parameters~$q$.
The values~$\mu$ may represent physical parameters or artificial parameters.

Thus the state variables or inner variables
$x: [0,t_{\rm end}] \times \mathcal{M} \rightarrow \R^n$
depend on time as well as the parameters.
The inputs $u: [0,t_{\rm end}] \rightarrow \R^{n_{\rm in}}$ are independent
of the parameters, whereas the outputs 
$y: [0,t_{\rm end}] \times \mathcal{M} \rightarrow \R^{n_{\rm out}}$
vary with respect to the parameters.
We consider a single output ($n_{\rm out} = 1$) without loss of generality.

We assume that the system~(\ref{ode-parameter}) is either
an ODE for all $\mu \in \mathcal{M}$ or a DAE for all $\mu \in \mathcal{M}$.
Let the system be asymptotically stable for each parameter
with respect to Definition~\ref{def:stable}.
Furthermore, initial value problems
\begin{equation} \label{ivp-parameter}
  x(0,\mu) = x_0(\mu)
\end{equation}
are considered including a function $x_0 : \mathcal{M} \rightarrow \R^n$.
The initial values may be independent of the parameters.
The initial values have to be consistent in the case of DAEs. 

We suppose that the parameters in the linear dynamical
system~(\ref{ode-parameter}) are affected by uncertainties.
A common approach is to substitute the parameters by
independent random variables
$\mu : \Omega \rightarrow \mathcal{M}$
on a probability space $(\Omega,\mathcal{F},P)$
with event space~$\Omega$, sigma-algebra~$\mathcal{F}$
and probability measure~$P$,
see~\cite{sullivan,xiu-book}. 
We apply traditional probability distributions like
uniform, Gaussian, beta, etc.
Hence a joint probability density function $\rho : \mathcal{M} \rightarrow \R$
is available.
This approach yields a stochastic model.

A measurable function $f: \mathcal{M} \rightarrow \R$ depending on the
random variables exhibits the expected value
\begin{equation} \label{expected-value}
  \mathbb{E} [ f ] = \displaystyle
  \int_{\Omega} f(\mu(\omega)) \; {\rm d}P(\omega) =
  \int_{\mathcal{M}} f(\mu) \rho(\mu) \; {\rm d}\mu
\end{equation}
provided that the integral is finite.
The expected value~(\ref{expected-value}) implies the inner product
\begin{equation} \label{inner-product}
  \langle f , g \rangle =
  \int_{\mathcal{M}} f(\mu) g(\mu) \rho(\mu) \; {\rm d}\mu
\end{equation}
for two functions in the Hilbert space
\begin{equation} \label{ltwo}
  \ltwo =
  \left\{ f : \mathcal{M} \rightarrow \R \; : \;
  f \; {\rm measurable}, \; \mathbb{E} [ f^2 ] < \infty \right\} .
\end{equation}
The associated norm is $\| f \|_{\ltwo} = \sqrt{\langle f , f \rangle}$
as usual.

\subsection{Polynomial chaos expansions}
In most cases, a complete orthogonal basis $(\Phi_i)_{i \in \N}$
of polynomials $\Phi_i : \mathcal{M} \rightarrow \R$ exists.
These multivariate polynomials are the products
\begin{equation} \label{basis-polynomials}
  \Phi_i(\mu) = \Psi_{i_1}^{(1)} (\mu_1) \cdot
  \Psi_{i_2}^{(2)} (\mu_2) \cdot \ldots \cdot
  \Psi_{i_{q-1}}^{(q-1)} (\mu_{q-1}) \cdot \Psi_{i_q}^{(q)} (\mu_{q}) ,
\end{equation}
where $(\Psi_{\ell}^{(j)})_{\ell \in \N_0}$ is the family of
univariate orthogonal polynomials
with respect to the $j$th random variable.
The degree of $\Psi_{\ell}^{(j)}$ is exactly~$\ell \ge 0$. 
Let this basis also be normalized.
Each traditional probability distribution implies its own family of
orthogonal basis polynomials, see~\cite{xiu-book}.
The basis is complete in the case of uniform, beta, and
Gaussian distribution, for example.
Yet the polynomials do not span the complete Hilbert space~(\ref{ltwo})
in the case of a log-normal distribution, see~\cite{ernst-etal}.

We assume that the QoI of the system~(\ref{ode-parameter}) is in the
space~(\ref{ltwo}) for each time point.
If the parameter domain~$\mathcal{M}$ is compact, then the
continuity of the QoI is sufficient for belonging to~(\ref{ltwo})
pointwise in time.
If the parameter domain~$\mathcal{M}$ is unbounded, then
integrability conditions have to be satisfied with respect to
the probability distribution. 
Consequently, the QoI can be expanded into the series
\begin{equation} \label{pce}
  y(t,\mu) = \sum_{i=1}^{\infty} w_i(t) \Phi_i(\mu) 
\end{equation}
with coefficient functions $w_i : [0,t_{\rm end}] \rightarrow \R$,
which is called a (generalized) polynomial chaos expansion (PCE).
The series~(\ref{pce}) converges in the norm of~(\ref{ltwo}) pointwise
in time.

Likewise, the state variables exhibit the PCE
\begin{equation} \label{pce-state}
  x(t,\mu) = \sum_{i=1}^{\infty} v_i(t) \Phi_i(\mu)
\end{equation}
with coefficient functions $v_i : [0,t_{\rm end}] \rightarrow \R^n$,
provided that each state variable is in the Hilbert space~(\ref{ltwo}).

We assume that the first basis polynomial is the unique constant
polynomial $\Phi_1 \equiv 1$.
The orthonormality $\langle \Phi_i , \Phi_j \rangle = \delta_{ij}$
of the basis functions imply the formulas
$$ \mathbb{E} [y(t,\cdot)] = w_1(t)
  \qquad \mbox{and} \qquad
  {\rm Var} [y(t,\cdot)] = \sum_{i=2}^{\infty} w_i(t)^2 $$
for the expected value and the variance of the QoI in each time point.

\subsection{Stochastic Galerkin method}
\begin{figure}
  \begin{center}  
  \includegraphics[width=7cm]{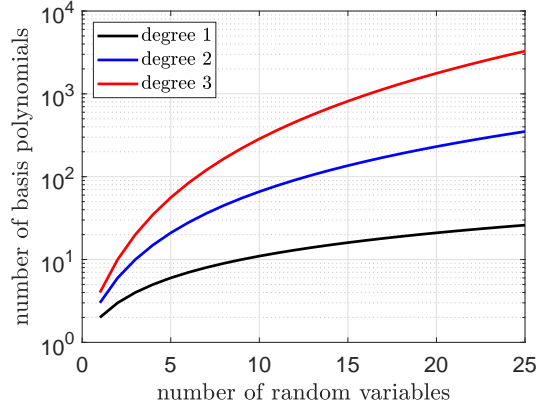}
  \end{center}
  \caption{Number of basis polynomials in dependence on
    number of random variables and total polynomial degree
    in semi-logarithmic scale.}
\label{fig:number-basis}
\end{figure}

We truncate the series~(\ref{pce}),(\ref{pce-state}) to a finite sum.
Typically, all basis polynomials up to some total degree~$d$ are included.
There is a bijective mapping between the integers~$i$
and the multi-indices $i_1,\ldots,i_q$
with the number~$q$ of random parameters.
In view of~(\ref{basis-polynomials}), we obtain the finite index set
$\{ i \in \N \, : \,  i_1 + i_2 + \cdots + i_{q-1} + i_q \le d \}$.
The cardinality of this index set is
$$ m = \frac{(d+q)!}{d!q!} . $$
If the number of random variables is large, then the number of basis
polynomials becomes huge even for moderate degrees, say $d = 3$.
Figure~\ref{fig:number-basis} illustrates the growth of the
number of basis polynomials.

Without loss of generality, the truncated series read as
\begin{equation} \label{truncated-pce}
  y^{(m)}(t,\mu) = \displaystyle{\sum_{i=1}^m} w_i(t) \Phi_i(\mu)
  \quad \mbox{and} \quad 
  x^{(m)}(t,\mu) = \displaystyle{\sum_{i=1}^m} v_i(t) \Phi_i(\mu) .
\end{equation}
Inserting the approximations~(\ref{truncated-pce}) into the
linear dynamical system~(\ref{ode-parameter}) yields a residual.
The Galerkin approach requires that the residual is orthogonal
to the space spanned by the basis polynomials $\Phi_1,\ldots,\Phi_m$
with respect to the inner product~(\ref{inner-product}). 
Basic calculations produce a larger coupled linear dynamical system
\begin{align} \label{galerkin}
  \begin{split}
  \hat{E} \dot{\hat{v}}(t) & = \hat{A} \hat{v}(t) + \hat{B} u(t) \\[1ex]
  \hat{w}(t) & = \hat{C} \hat{v}(t)
  \end{split}
\end{align} 
for $\hat{v} = (\hat{v}_1^\top , \ldots , \hat{v}_m^\top)^\top$ and
$\hat{w} = (\hat{w}_1,\ldots,\hat{w}_m)^\top$.
Hence we obtain a linear dynamical system of dimension $\hat{n} = mn$
with $m$~outputs.
The matrices exhibit the sizes $\hat{A},\hat{E} \in \R^{\hat{n} \times \hat{n}}$,
$\hat{B} \in \R^{\hat{n} \times n_{\rm in}}$, and $\hat{C} \in \R^{m \times \hat{n}}$.
To define the matrices, we introduce the auxiliary arrays
\begin{equation} \label{auxiliary-array}
  {S}(\mu) = (\Phi_i(\mu)\Phi_j(\mu))_{i,j=1,\ldots,m}
  \qquad \mbox{and} \qquad
  {s}(\mu) = (\Phi_i(\mu))_{i=1,\ldots,m} .
\end{equation}
It follows that
\begin{equation} \label{galerkin-matrices}
\hat{A} = \mathbb{E} [{S} \otimes A] , \quad
\hat{B} = \mathbb{E} [{s} \otimes B] , \quad
\hat{C} = \mathbb{E} [{S} \otimes C], \quad
\hat{E} = \mathbb{E} [{S} \otimes E]
\end{equation}
using Kronecker products.
Therein, the probabilistic integration~(\ref{expected-value}) operates
separately in each component of the matrices.
Often it holds that ${\rm rank}(\hat{C}) = m$.
More details on the stochastic Galerkin method for linear dynamical systems
are given in~\cite{pulch14,pulch-maten}.

The mass matrix~$\hat{E}$ may be singular even though all matrices~$E$ are
non-singular due to properties of the spectrum, see~\cite{sonday}.
However, this loss of invertibility hardly occurs.
Yet the mass matrix~$\hat{E}$ may become ill-conditioned.
We assume that $\hat{E}$ is always non-singular in the case of
ODEs~(\ref{ode-parameter}).
Likewise, the stochastic Galerkin system~(\ref{galerkin}) may be unstable
even though the systems~(\ref{ode-parameter}) are asymptotically stable
for all parameters.
Academic examples are given in~\cite{pulch-augustin}.
Again this loss of stability is hardly observed in practice. 
Furthermore, the passivity of stochastic Galerkin systems was
investigated for models of electric circuits in~\cite{manfredi}.

Initial conditions $\hat{v}(0)=\hat{v}_0$ are derived from the
initial values~(\ref{ivp-parameter}) of the
original dynamical system~(\ref{ode-parameter}) by an own truncated PCE.
If the initial values~(\ref{ivp-parameter}) are identical to zero,
then the choice $\hat{v}(0)=0$ is obvious.
The approximation of the QoI reads as 
\begin{equation} \label{appr-qoi}
  \hat{y}^{(m)}(t,\mu) =
  \displaystyle{\sum_{i=1}^m} \hat{w}_i(t) \Phi_i(\mu) ,
\end{equation}
where the outputs $\hat{w}_1,\ldots,\hat{w}_m$
of the stochastic Galerkin system~(\ref{galerkin}) yield the coefficients.

If the entries of the matrices are polynomials depending on~$\mu$
in the system~(\ref{ode-parameter}), then the
matrices~(\ref{galerkin-matrices}) of the stochastic Galerkin system can
be calculated analytically.
Consequently, we obtain the matrices exactly
(except for round-off errors) independent of the number~$q$ of
random variables.
In contrast, stochastic collocation techniques,
which are non-intrusive methods, cf.~\cite{sullivan,xiu-book},
induce a quadrature error or sampling error.
This error typically grows for fixed numbers of collocation points and
increasing dimensions~$q$.
Although the stochastic Galerkin approach represents an intrusive method,
the effort of coding the algorithms is not extensive,
because just constant matrices have to be specified for
linear time-invariant systems.

We prove a property, which will be used later.
\begin{lemma} \label{lemma:definite}
  If the matrices $E(\mu)$ are symmetric and positive definite for
  almost all $\mu \in \mathcal{M}$,
  then the stochastic Galerkin projection~$\hat{E}$ is also symmetric 
  and positive definite.
\end{lemma}

\underline{Proof:}

The Galerkin-projected matrix consists of the blocks 
$\hat{E}_{ij} = \mathbb{E}[\Phi_i \Phi_j E]$ for $i,j=1,\ldots,m$,
see~(\ref{galerkin-matrices}).
Hence the symmetry is obvious.
Let $z = (z_1^\top,\ldots,z_m^\top)^\top \in \R^{mn}$.
We obtain
\begin{align*}
  \begin{split}
  z^\top \hat{E} z & = \displaystyle
  \sum_{i,j=1}^m z_i^\top \hat{E}_{ij} z_j
  =
  \sum_{i,j=1}^m z_i^\top \mathbb{E}[\Phi_i \Phi_j E] z_j
  =
  \mathbb{E} \left[ \sum_{i,j=1}^m z_i^\top E z_j \Phi_i \Phi_j\right] \\
  & = \displaystyle
  \mathbb{E} \left[ \left( \sum_{i=1}^m z_i \Phi_i \right)^\top
    E \left( \sum_{j=1}^m z_j \Phi_j \right) \right]
  \;\; \ge \;\; 0 ,
  \end{split}
\end{align*} 
because the integrand is almost everywhere non-negative
in the probabilistic integration~(\ref{expected-value}).
The basis functions $(\Phi_i)_{i \in \N}$ are linearly independent.
Thus $z \neq 0$ implies that the above finite sum is non-zero on a
subset $\mathcal{U} \subset \mathcal{M}$ satisfying
$P(\{ \omega : \mu(\omega) \in \mathcal{U}\}) > 0$
for the probability measure~$P$.
It follows that $z^\top \hat{E} z > 0$ for $z \neq 0$.
\hfill $\Box$

\medskip

Likewise, this relation applies to the case of negative definite matrices.


\section{Stability preservation}
\label{sec:stability}
\begin{figure}
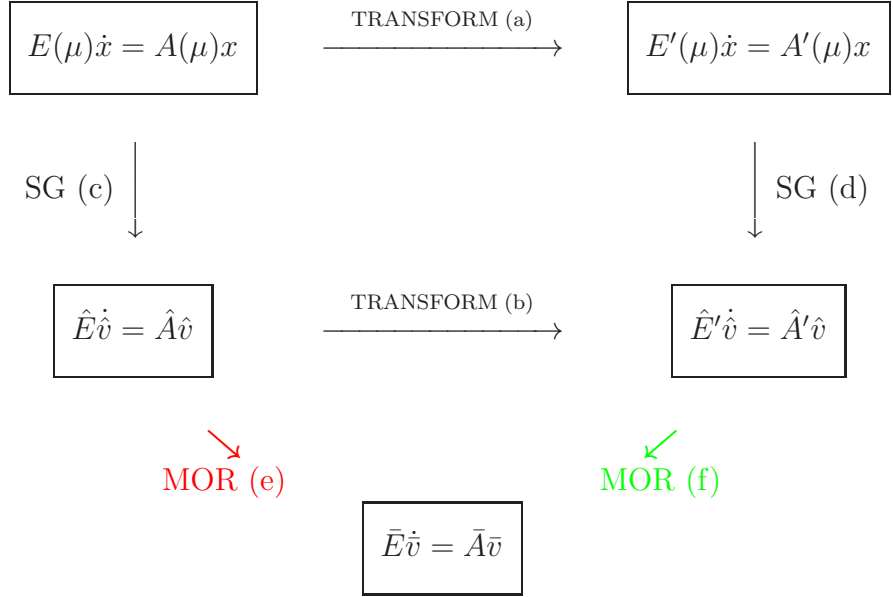

\begin{center}
  \begin{tabular}{ccccc}
    \fbox{$\; E(\mu) \dot{x} = A(\mu) x$ \rule[-4mm]{0mm}{10mm}} & &
    $\stackrel{\rm TRANSFORM \; (a)}{\xrightarrow{\hspace*{3cm}}}$ & &
    \fbox{$\; {E}'(\mu) \dot{x} = {A}'(\mu) x$ \rule[-4mm]{0mm}{10mm}} \\[1cm]
    SG (c) $\boldsymbol{\Bigg\downarrow}$ \hspace{1cm} \mbox{} & & & &
    \mbox{} \hspace{1cm} $\boldsymbol{\Bigg\downarrow}$ SG (d) \\[1.0cm]
    \fbox{$\; \hat{E} \dot{\hat{v}} = \hat{A} \hat{v}$
      \rule[-4mm]{0mm}{10mm}} & &
    $\stackrel{\rm TRANSFORM \; (b)}{\xrightarrow{\hspace*{3cm}}}$ & &
    \fbox{$\; \hat{E}' \dot{\hat{v}} = \hat{A}' \hat{v}$
      \rule[-4mm]{0mm}{10mm}} \\[1.0cm]
    & \mbox{} \hspace{-2cm} \color{red} $\boldsymbol{\searrow}$ & &
    \color{green} $\boldsymbol{\swarrow}$ \color{red} \hspace{-2cm} \mbox{}
    & \\
    & \mbox{} \hspace{-2cm} \color{red} MOR (e) & &
    \color{green} MOR (f) \color{black}  \hspace{-2cm} \mbox{} & \\
    & &
    \fbox{$\; \bar{E} \dot{\bar{v}} = \bar{A} \bar{v}$ \rule[-4mm]{0mm}{10mm}}
    & & 
  \end{tabular}
\end{center}
\caption{Flowchart for transformations,
  stochastic Galerkin (SG) projections and model order reduction (MOR).
  Inputs and outputs are not shown. \label{fig:flowchart}}
\end{figure}

We investigate three strategies to preserve the asymptotic stability
in an MOR of a stochastic Galerkin system.
Figure~\ref{fig:flowchart} illustrates different possibilities.

\subsection{Transformation of stochastic Galerkin system}
\label{sec:trafo-gal}
This approach follows the steps (c),(b),(f) in Figure~\ref{fig:flowchart}.
If the linear dynamical systems~(\ref{ode-parameter}) are asymptotically
stable for all parameters, then the stochastic Galerkin
system~(\ref{galerkin}) is usually asymptotically stable.
However, if the system~(\ref{galerkin}) is not dissipative, then
stability may be lost in some MOR methods by step~(e). 
In this case, we can transform the system to a dissipative form
as demonstrated in Section~\ref{sec:transformation}
and Section~\ref{sec:solve-lyapunov}.
Consequently, the reduced system is stable due to
Theorem~\ref{thm:dissipative}.
Remark that we do not have to calculate the transformed matrices
of the high-dimensional system explicitly.
Just an appropriate projection matrix has to be determined
via~(\ref{projection-w}).

In this approach, the critical part is the solution of the
high-dimensional Lyapunov equation~(\ref{lyapunov}).
A direct method of linear algebra would require
$\mathcal{O}(m^3n^3)$ operations.
Thus we are restricted to approximate methods or iteration schemes.
The possibilities are listed in Section~\ref{sec:solve-lyapunov}.

\subsection{Transformation of parameter-dependent system}
\label{sec:trafo-org}
Now the succession (a),(d),(f) is performed with respect
to Figure~\ref{fig:flowchart}. 
Consequently, we transform the linear dynamical systems~(\ref{ode-parameter})
first.

\subsubsection{Transformation}
The stochastic Galerkin system is not always asymptotically stable,
even if all parameter-dependent systems~(\ref{ode-parameter}) are
asymptotically stable. 
However, we obtain a positive result in the case of dissipativity.

\newpage

\begin{theorem} \label{thm:gal-diss}
  If the linear dynamical systems~(\ref{ode-parameter}) are dissipative
  for almost all $\mu \in \mathcal{M}$,
  then the stochastic Galerkin system~(\ref{galerkin}) is also dissipative.
  Consequently, a Galerkin-type reduction of~(\ref{galerkin}) yields
  an asymptotically stable system.
\end{theorem}

\underline{Proof:}

Lemma~\ref{lemma:definite} implies that the mass matrix~$\hat{E}$ is
symmetric and positive definite.
The matrix $\hat{A} + \hat{A}^\top$ is the Galerkin projection
of $A(\cdot) + A(\cdot)^\top$ due to
$$ \hat{A} + \hat{A}^\top =
\mathbb{E} [{S} \otimes A] +
\mathbb{E} [{S} \otimes A^\top] =
\mathbb{E} [{S} \otimes A + {S} \otimes A^\top] =
\mathbb{E} [{S} \otimes (A+A^\top)] , $$
see~(\ref{galerkin-matrices}).
Since $A(\mu)+A(\mu)^\top$ is negative definite for almost
all $\mu \in \mathcal{M}$,
Lemma~\ref{lemma:definite} shows that $\hat{A} + \hat{A}^\top$ is
negative definite.
Hence both conditions of Definition~\ref{def:dissipative} are satisfied.
Theorem~\ref{thm:dissipative} guarantees the stability of a reduced system.
\hfill $\Box$

\medskip

If the original systems~(\ref{ode-parameter}) are not dissipative for almost
all $\mu \in \mathcal{M}$, then we transform the systems appropriately.
The following transformation has to be applied for almost all~$\mu$
simultaneously (even if some system is already dissipative)
to preserve continuity and smoothness of the functions.

We obtain the parameter-dependent Lyapunov equation
\begin{equation} \label{lyap-par}
    A(\mu)^\top M(\mu) E(\mu) + E(\mu)^\top M(\mu) A(\mu) + F = 0
\end{equation}
for $\mu \in \mathcal{M}$.
Although the matrix~$F$ may depend on the parameters, 
we choose a constant matrix,
because no parameter-aware strategy with benefits is known yet.
The Lyapunov equations~(\ref{lyap-par}) yield a unique family~$M(\mu)$
of symmetric positive definite matrices.

Furthermore, the stochastic Galerkin system~(\ref{galerkin})
obtained by step~(c) is not equivalent to a stochastic Galerkin system
obtained by steps~(a),(d).
On the one hand, the stochastic Galerkin method is invariant with respect to
transformations by constant non-singular matrices.
On the other hand, we use the parameter-dependent matrix~$E(\mu)^\top M(\mu)$
in the transformation~(a) to the dissipative form~(\ref{ode-diss}).

\subsubsection{Polynomial system matrices}
Now we assume that the matrices $A(\mu),B(\mu),E(\mu)$
of the system~(\ref{ode-parameter})
involve only polynomials in the variable~$\mu$,
which is often given in practice.
Thus the matrices $\hat{A},\hat{B},\hat{E}$ of the
stochastic Galerkin system~(\ref{galerkin})
can be calculated analytically in the case of
traditional probability distributions.
We do not consider the output matrix $C(\mu)$,
because it is not transformed. 
However, the matrix-valued function $M(\mu)$ satisfying~(\ref{lyap-par})
consists of rational functions in the variable~$\mu$.
In the case of low dimensions~$n$,
the function $M(\mu)$ can be calculated explicitly by a
computer algebra software.
Alternatively, the solution $M(\mu)$ can be evaluated for
a finite set of parameters~$\mu$.

We require the Galerkin projection of the transformed matrices
$$ A'(\mu) = E(\mu)^\top M(\mu) A(\mu) , \;
B'(\mu) = E(\mu)^\top M(\mu) B(\mu) , \;
E'(\mu) = E(\mu)^\top M(\mu) E(\mu) $$
in the dissipative system~(\ref{ode-diss}).
The entries of these transformed matrices are rational functions of~$\mu$.
A quadrature rule yields numerical approximations
$\tilde{A},\tilde{B},\tilde{E}$ of the matrices
$\hat{A}',\hat{B}',\hat{C}'$.
Let $\{ \mu_1,\ldots,\mu_k \} \subset \mathcal{M}$ be the nodes
and $\{ \gamma_1,\ldots,\gamma_k \} \subset \R$ be the weights.
The approximation of the blocks in~$\hat{A}'$
reads as
\begin{equation} \label{quadrature}
  \hat{A}_{ij}' = \mathbb{E}[\Phi_i \Phi_j A'] \approx
  \sum_{\ell=1}^k \gamma_{\ell} \Phi_i(\mu_{\ell}) \Phi_j(\mu_{\ell})
  E(\mu_{\ell})^\top M(\mu_{\ell}) A(\mu_{\ell})
\end{equation}
for $i,j=1,\ldots,m$.
Using the auxiliary array~(\ref{auxiliary-array}), the approximation becomes
\begin{equation} \label{quadrature-matrix}
  \tilde{A} =
  \sum_{\ell=1}^k \gamma_{\ell} \; S(\mu_{\ell}) \otimes {A}'(\mu_{\ell}) .
\end{equation}
Likewise, the quadrature yields $\tilde{B}$ and $\tilde{E}$.
The computational effort is dominated by determining 
$k$~numerical solutions of the Lyapunov equations~(\ref{lyap-par}).

\begin{lemma} \label{lemma:quadrature}
  If all weights are positive,
  then a quadrature of kind~(\ref{quadrature}) yields
  approximations~$\tilde{A},\tilde{E}$,
  where $\tilde{E}$ is symmetric positive semi-definite and
  $\tilde{A} + \tilde{A}^{\top}$ is negative semi-definite.
\end{lemma}

\underline{Proof:}

The symmetry of $\tilde{E}$ is obvious.
Let $z \in \R^{mn} \backslash \{ 0 \}$
with $z = (z_1^\top,\ldots,z_m^\top)^\top$ and
$$ v_{\ell} =  \sum_{j=1}^m z_j \Phi_j(\mu_{\ell})
\qquad \mbox{for} \;\; \ell = 1,\ldots,k . $$
We obtain using the notation~(\ref{quadrature-matrix})
\begin{align*}
  \begin{split}
    z^\top \tilde{E} z & = \sum_{\ell=1}^k \gamma_{\ell}
    z^\top ( S(\mu_{\ell}) \otimes {E}'(\mu_{\ell}) ) z 
    = \sum_{\ell=1}^k \gamma_{\ell} \sum_{i,j=1}^m
    \Phi_i(\mu_{\ell}) \Phi_j(\mu_{\ell}) z_i^\top {E}'(\mu_{\ell}) z_j \\
    & = \sum_{\ell=1}^k \gamma_{\ell} v_{\ell}^\top {E}'(\mu_{\ell}) v_{\ell} . \\
  \end{split}
\end{align*}
Since $E'(\mu_{\ell})$ is positive definite and the weights~$\gamma_{\ell}$
are positive for each~$\ell$,
all terms of the sum are non-negative.
Hence $\tilde{E}$ is positive semi-definite.
The negative semi-definiteness of $\tilde{A} + \tilde{A}^{\top}$
is concluded by the same treatment.
\hfill $\Box$

\medskip

In addition, it is very likely that one or more terms are positive
in the above sum.
Thus we assume that these matrices are definite.
It follows that the approximate Galerkin system has the
dissipativity condition of Definition~\ref{def:dissipative}.

Concerning the positivity of the weights, we address three classes
of multivariate quadrature schemes or sampling methods,
which can all be found in~\cite{sullivan}:
\begin{itemize}
\item[i)] {\em tensor-product formulas}:
  Univariate quadrature rules often exhibit exclusively positive weights
  like Gaussian quadrature, for example.
  The tensor-product formulas inherit the positivity,
  because their weights are the products of the weights
  in the univariate case.
\item[ii)] {\em sparse grid quadrature}:
  Often both positive and negative weights occur.
  Sparse grids with purely positive weights typically require more nodes
  for the same accuracy.
\item[iii)] {\em (quasi) Monte-Carlo methods}:
  In these sampling techniques, 
  the weights are $\gamma_{\ell}=\frac{1}{k}$ for all $\ell=1,\ldots,k$
  and thus positive.
\end{itemize}

\subsubsection{Properties}
\label{sec:par-trafo-properties}
The strategy of this section looks appealing
in the case of low dimensions~$n$,
because the Lyapunov equations can be solved cheap.
However, the number of random parameters is typically large in this case,
because otherwise the stochastic Galerkin system is not high-dimensional
and an MOR is obsolete.
A large number of random variables implies a quadrature in a
high-dimensional space.
Using a computationally feasible number of nodes,
the quadrature error may still be too large such that the exact solution
of the approximate Galerkin system yields poor approximations~(\ref{appr-qoi})
of the random-dependent QoI.
Therefore, the above approach is critical.

Furthermore, there is a major loss of sparsity in this technique.
Let the entries in the matrices be polynomials of low degrees
depending on~$\mu$ in the system~(\ref{ode-parameter}).
Even if the matrices are dense, then the stochastic Galerkin
system~(\ref{galerkin}) exhibits sparse matrices due to the
orthogonality of the basis polynomials.
However, the dense matrices of the transformed systems~(\ref{ode-diss})
are rational functions depending on~$\mu$.
Inner products~(\ref{inner-product}) of their entries and
orthogonal polynomials are non-zero and thus the matrices
of~(\ref{galerkin}) become dense in the Galerkin projection.
In contrast to the approach of Section~\ref{sec:trafo-gal},
we have to calculate the matrices of the alternative Galerkin system
explicitly to determine a projection matrix~$V$ of the MOR.
 
\subsection{Transformation using reference parameter}
\label{sec:ref-parameter}
Again the steps (c),(b),(f) are considered in Figure~\ref{fig:flowchart},
where the transformation is done different from
Section~\ref{sec:trafo-gal}.
We derive an additional technique to decrease the computational effort
and to omit quadrature errors.
A single reference parameter~$\mu^* \in \mathcal{M}$ is selected
like the mean value $\mu^* = \mathbb{E}[\mu]$, for example.
We directly solve the Lyapunov equation~(\ref{lyap-par})
only for~$\mu^*$ using some matrix~$F$.
The solution $M^* = M(\mu^*) \in \R^{n \times n}$ is symmetric and
positive definite.
We define the larger transformation matrix
\begin{equation} \label{large-M}
  \hat{M} = I_m \otimes M^* \in \R^{mn \times mn}
\end{equation}
using the identity matrix $I_m \in \R^{m \times m}$ and the Kronecker product.
Obviously, the matrix~(\ref{large-M}) is symmetric and
positive definite again.
We employ this matrix to transform the stochastic Galerkin
system~(\ref{galerkin}) into the form~(\ref{ode-diss}).
Since the matrix~(\ref{large-M}) is block-diagonal,
matrix-matrix multiplications with $\hat{M}$ are cheap.
We do not need to compute the transformed system matrices.
Alternatively, we directly compute the projection matrix~(\ref{projection-w}),
where a projection matrix~$V$ is determined by the
original Galerkin system.

Any random variables can be decomposed into the
form $\mu(\omega) = \mu^* + \Delta \mu(\omega)$.
We consider the random variables
\begin{equation} \label{random-variables}
  \mu_\theta(\omega) = \mu^* + \theta \Delta \mu(\omega)
\end{equation}
with the real parameter $\theta \in [0,1]$.

\begin{theorem} \label{thm:limitcase}
  Let the random variables in a linear system of ODEs~(\ref{ode-parameter})
  be of the form (\ref{random-variables}).
  There is a positive constant $\theta_0 \in (0,1]$ such that
    the transformation of a stochastic Galerkin system~(\ref{galerkin})
    using the matrix~(\ref{large-M}) yields a dissipative system
    of type~(\ref{ode-diss}) for
    all $\theta \in [0,\theta_0]$. 
\end{theorem}

\underline{Proof:}

The transformed mass matrix $\hat{E}^\top \hat{M} \hat{E}$ is 
symmetric and positive definite for all~$\theta$,
since just the non-singularity of the original mass matrix~$\hat{E}$
is required.
In the limit, we obtain
\begin{equation} \label{trafo-limit}
\lim_{\theta \rightarrow 0} \hat{E}^\top \hat{M} \hat{A}
= I_m \otimes (E(\mu^*)^\top M^* A(\mu^*)) 
\end{equation}
due to the orthogonality of the basis polynomials.
Hence the matrix becomes block-diagonal with identical blocks.
The factor~$M^*$ satisfies the Lyapunov equation~(\ref{lyap-par})
for the chosen positive definite matrix~$F$.
The symmetric part of the matrix~(\ref{trafo-limit}) is negative definite, 
because it holds that
$$ \begin{array}{clcl}
  & I_m \otimes (E(\mu^*)^\top M^* A(\mu^*)) +
  \big( I_m \otimes (E(\mu^*)^\top M^* A(\mu^*)) \big)^\top & & \\[1ex]
  = &
  I_m \otimes (E(\mu^*)^\top M^* A(\mu^*)) +
  I_m \otimes (A(\mu^*)^\top M^* E(\mu^*))
  & = & - I_m \otimes F . \\
\end{array} $$
The two conditions of Definition~\ref{def:dissipative} are satisfied
and thus the system is dissipative in the limit.
Since the limit~(\ref{trafo-limit}) is reached continuously
with respect to the parameter~$\theta$,
it follows that a sufficiently small perturbation of $\theta=0$ still
yields matrices with the required properties.
\hfill $\Box$

\medskip

Theorem~\ref{thm:limitcase} shows that the stochastic Galerkin system
becomes dissipative for all sufficiently small~$\theta > 0$.
However, this implication does not guarantee that the system is dissipative
for $\theta=1$, which reproduces our desired choice of random variables
in~(\ref{random-variables}).
Nevertheless, the computational effort is low such that it is worth
to try this approach.
Even if this transformed stochastic Galerkin system is not dissipative,
a loss of stability may happen less often in an MOR.

\subsection{Nonlinear dynamical systems}
We outline a stabilization concept for nonlinear dynamical systems.
Let an autonomous system of ODEs be given in the form
\begin{align} \label{nonlinear-ode}
  \begin{split}
  E(\mu) \dot{x}(t,\mu) & = f(x(t,\mu),\mu) \\[1ex]
  y(t,\mu) & = g(x(t,\mu),\mu)
  \end{split}
\end{align}
with a parameter-dependent mass matrix~$E \in \R^{n \times n}$
and a nonlinear smooth right-hand side
$f : \R^n \times \mathcal{M} \rightarrow \R^n$.
The QoI~$y$ is defined by the linear or nonlinear function
$g : \R^n \times \mathcal{M} \rightarrow \R^{n_{\rm out}}$.
Let $n_{\rm out} = 1$.
We assume that a family of asymptotically stable stationary solutions
$x^* : \mathcal{M} \rightarrow \R^n$ exists, i.e.,
$f(x^*(\mu),\mu) = 0$ for all $\mu \in \mathcal{M}$.
The definition of stable stationary solutions can be found
in~\cite[p.~22]{seydel}.
As in~\cite{pulch-augustin},
the system~(\ref{nonlinear-ode}) is replaced by the equivalent system
\begin{align} \label{nonlinear-ode-trafo}
  \begin{split}
    E(\mu) \dot{\tilde{x}}(t,\mu) & =
    \tilde{f}(\tilde{x}(t,\mu),\mu) = 
    f(\tilde{x}(t,\mu)+x^*(\mu),\mu) \\[1ex]
    y(t,\mu) & =
    \tilde{g}(\tilde{x}(t,\mu),\mu) = 
    g(\tilde{x}(t,\mu)+x^*(\mu),\mu) , 
  \end{split}
\end{align}
which exhibits the constant asymptotically stable stationary solution
$\tilde{x}^* = 0$.
Let $\tilde{J}(\mu) \in \R^{n \times n}$ be the Jacobian matrix of $\tilde{f}$
evaluated at $\tilde{x}^*=0$ and~$\mu \in \mathcal{M}$.
The stochastic Galerkin method yields a larger nonlinear system of ODEs
\begin{align} \label{galerkin-nonlinear}
  \begin{split}
  \hat{E} \dot{\hat{v}}(t) & = \hat{f}(\hat{v}(t)) \\[1ex]
  \hat{w}(t) & = \hat{g}(\hat{v}(t)) 
  \end{split}
\end{align}
including the mass matrix~$\hat{E}$ from~(\ref{galerkin-matrices}).
The definition of the nonlinear function
$\hat{f} : \R^{mn} \rightarrow \R^{mn}$ is given
in~\cite{pulch19,pulch-augustin}.
This function inherits the smoothness of~$f$. 
Now $\hat{v}^*=0$ is a stationary solution of~(\ref{galerkin-nonlinear}).
Although this stationary solution may be unstable,
this loss of stability hardly occurs in practice.
Thus we assume that the stationary solution is asymptotically stable again.
Let $\hat{J} \in \R^{mn \times mn}$ be the Jacobian matrix of~$\hat{f}$
evaluated at~$\hat{v}^*=0$.
It follows that the spectral abscissa of the matrix pencil $(\hat{E},\hat{J})$
is negative, see Definition~\ref{def:abscissa}.

A Galekin-type projection-based MOR yields a small dynamical system
\begin{align} \label{rom-nonlinear}
  \begin{split}
    \bar{E} \dot{\bar{v}}(t) & = \bar{f}(\bar{v}(t)) =
    V^\top \hat{f}(V \bar{v}(t)) \\[1ex]
  \bar{w}(t) & = \bar{g}(\bar{v}(t)) = \hat{g}(V \bar{v}(t)) 
  \end{split}
\end{align}
with a projection matrix~$V$ and the mass matrix $\bar{E} = V^\top \hat{E} V$. 
The reduced system owns the stationary solution $\bar{v}^* = 0$,
which may be unstable.

Stability-preserving techniques can be derived.
We consider the method of Section~\ref{sec:trafo-gal}, 
where the high-dimensional Lyapunov equation~(\ref{lyapunov})
is solved including the mass matrix~$\hat{E}$ from~(\ref{galerkin-nonlinear})
and $A = \hat{J}$.
The Galerkin system~(\ref{galerkin-nonlinear}) is transformed and
reduced to~(\ref{rom-nonlinear}), while the asymptotic stability of the
stationary solution is guaranteed as shown in~\cite[p.~35]{pulch-naco}
for general nonlinear implicit systems of ODEs.
Concerning the stability-preserving approach of Section~\ref{sec:trafo-org},
parameter-dependent Lyapunov equations~(\ref{lyap-par}) are solved
including $E(\mu)$ from~(\ref{nonlinear-ode}) and
$A(\mu) = \tilde{J}(\mu)$ with the Jacobian matrix associated
to the stationary solution of~(\ref{nonlinear-ode-trafo}).
This technique was used for explicit ODEs in~\cite{pulch-augustin}.
The application of the stability-preserving method from
Section~\ref{sec:ref-parameter} is straightforward now.


\section{Differential-algebraic equations}
If the linear dynamical system~(\ref{ode}) features a singular
mass matrix~$E$, then a system of DAEs is given.
MOR methods are also available for DAEs, see~\cite{benner-stykel}.
Yet the Lyapunov equations~(\ref{lyapunov}) do not have a solution,
which is crucial in our stability-preserving technique.
Let $x \in \R^n$.
It follows that
$$ x^\top A^\top M (Ex) + (Ex)^\top M Ax = - x^\top F x . $$
Choosing an $x \neq 0$ in the kernel of the matrix~$E$ implies
$Ex=0$ and thus a contradiction to the definiteness of the
matrix~$F$ appears.
We require an alternative strategy now.

\subsection{Regularization}
We apply a regularization of an asymptotically stable DAE system~(\ref{ode}),
which was also used in~\cite{mueller}.
The asymptotic stability guarantees a regular matrix pencil
as specified by Definition~\ref{def:matrix-pencil}.
The regularized system matrices read as
\begin{align} \label{regularization}
  \begin{split}
    E_{\rm reg} & = E - \alpha A \\
    A_{\rm reg} & = A + \beta E
  \end{split}
\end{align}
introducing parameters $\alpha,\beta>0$.
The matrix~$E_{\rm reg}$ is non-singular for all $\alpha > 0$.
Choosing $\alpha = \beta^2$,
it follows that the linear dynamical system~(\ref{ode})
with $A_{\rm reg},E_{\rm reg}$ is asymptotically stable for
all sufficiently small parameters~$\beta$, see~\cite{mueller}. 
An advantage of this regularization technique is that
the sparsity pattern of $s E - A$ coincides with the
sparsity pattern of $s E_{\rm reg} - A_{\rm reg}$ for each $s \in \C$.
Hence no loss of sparsity happens in the relevant operations.

We also choose a small parameter~$\beta$ (and $\alpha = \beta^2$) to ensure
that the difference between the original DAE system and the regularized ODE
system is small.
Assuming that the DAE system together with the defined outputs
has a transfer function with finite $\htwo$-norm~(\ref{htwo-norm}), 
error bounds were derived with respect to this norm in~\cite{pulch-arxiv}.

\subsection{Stochastic Galerkin projection}
\label{sec:dae-gal}
The steps (a)-(d) refer to the flowchart in Figure~\ref{fig:flowchart}.
The following two approaches are equivalent
provided that the same parameters $\alpha,\beta$ are chosen
in a regularization~(\ref{regularization}):
\begin{itemize}
\item[i)]
  Regularize the parameter-dependent system~(\ref{ode-parameter})
  and then project the systems of ODEs
  in the stochastic Galerkin method. 
\item[ii)]
  Project the parameter-dependent DAE system~(\ref{ode-parameter})
  in the stochastic Galer\-kin method and then 
  regularize the Galerkin system~(\ref{galerkin}).
\end{itemize}
In both cases, the stochastic Galerkin projection~(c) yields
the same matrices $\hat{A}$ and $\hat{E}$ for the system~(\ref{galerkin}). 
The mass matrix~$\hat{E}$ is non-singular due to the
regularization. 

The transformations (a) and (b) 
to a dissipative representation are done only for a regularized system,
since a system of ODEs is required for each Lyapunov equation.
The transformation~(b) involves the matrices from~(i) or~(ii).
However, in the transformation~(a), we have to regularize
the systems of DAEs~(\ref{ode-parameter}) immediately.
Furthermore, the drawbacks of the succession~(a),(d),
see Section~\ref{sec:par-trafo-properties},
also apply in this case.

In the transformation~(b), the approximate solution of the
high-dimensional Lyapunov equations~(\ref{lyapunov}) may
become critical.
Small regularization parameters~$\alpha$ 
nearly coincide with the singular case.
This problem is less pronounced in direct methods to solve
the Lyapunov equations.

We can also employ the technique from Section~\ref{sec:ref-parameter}
in this context.
The linear dynamical system~(\ref{ode-parameter}) is regularized for a
reference parameter~$\mu^*$ and some $\alpha,\beta>0$.
Solving the Lyapunov equation~(\ref{lyap-par}) including~$\mu^*$ yields the
high-dimensional transformation matrix~(\ref{large-M}).
Now the regularized Galerkin system from strategy~(ii),
where the same values $\alpha,\beta$ are used,
is transformed based on~(\ref{large-M}).
Again the dissipativity of the transformed representation cannot be
guaranteed a priori.


\section{Illustrative examples}
\label{sec:example}
We investigate a system of ODEs as well as a system of DAEs
in this section.
We computed on a FUJITSU Esprimo P920
Intel(R) Core(TM) i5-4570 CPU with 3.20 GHz (4~cores)
and operation system Microsoft Windows~7.
The software package MATLAB~\cite{matlab2018} (version R2018a)
was used for all computations.

\subsection{Mass-spring-damper system}
\label{sec:msd}
\begin{figure}
  \begin{center}  
  \includegraphics[width=5cm]{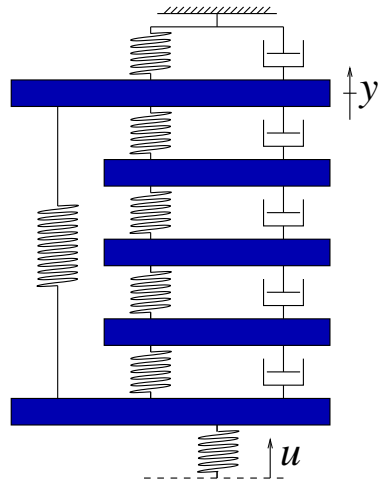}
  \end{center}
  \caption{Mass-spring-damper configuration.}
\label{fig:massspringdamper}
\end{figure}

\begin{figure}
  \begin{center}  
  \includegraphics[width=6.5cm]{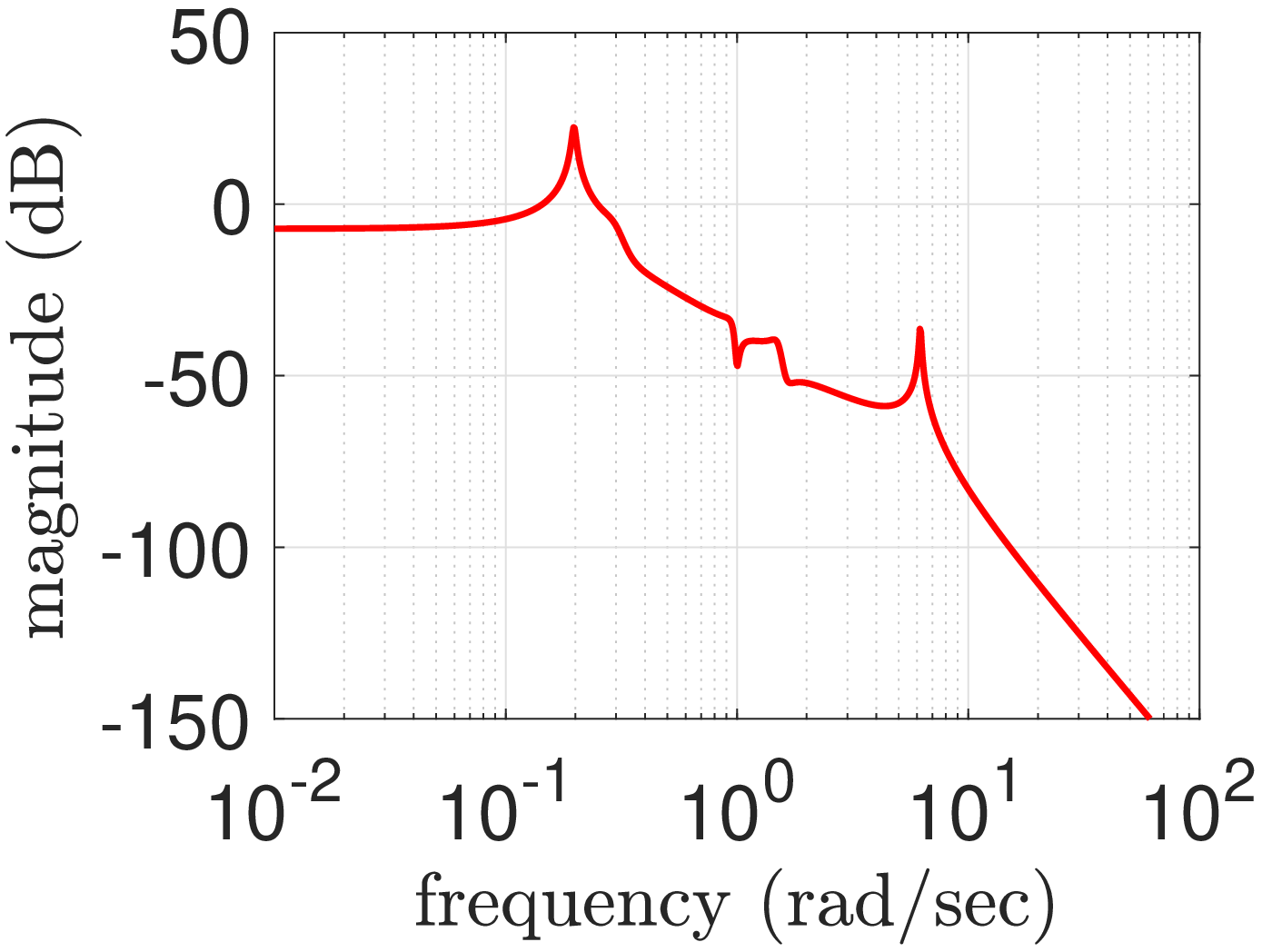}
  \hspace{5mm}
  \includegraphics[width=6.5cm]{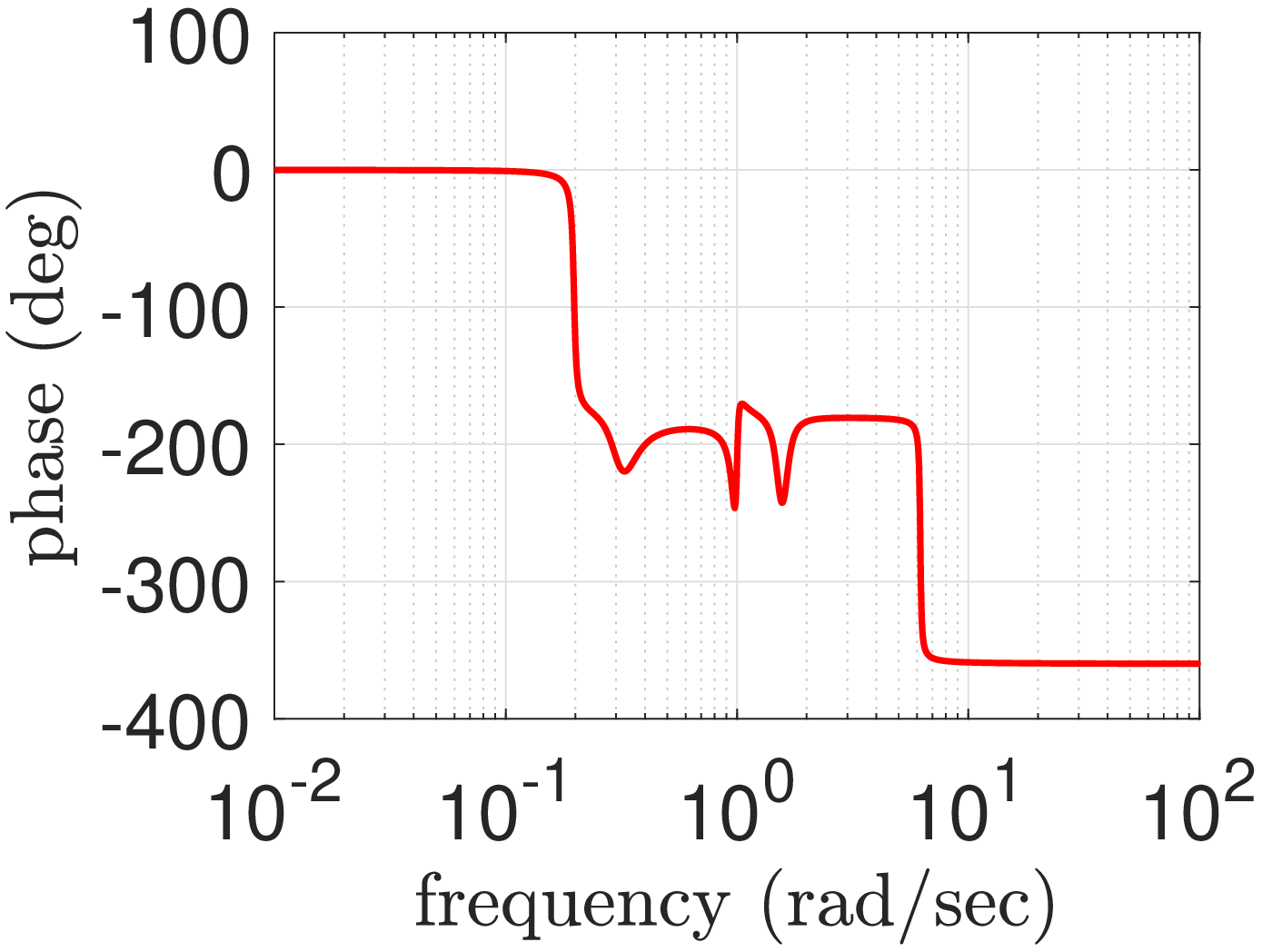}
  \end{center}
  \caption{Bode plot of mass-spring-damper system for
  a constant choice of the parameters.}
\label{fig:msd-bode}
\end{figure}

Figure~\ref{fig:massspringdamper} depicts a mechanical
configuration, which consists of 5~masses, 7~springs and 5~dampers.
The single input is the excitation at the bottom spring,
whereas the single output is the position of the top mass.
A mathematical modeling yields a linear dynamical system~(\ref{ode-parameter})
of $n=10$ ODEs of first order including $q=17$ parameters.
The system is asymptotically stable for all positive parameters.
The system matrices are affine-linear functions (polynomials of degree one)
of the parameters. 
The Bode plot of the system is shown for a specific choice of the
parameters in Figure~\ref{fig:msd-bode}.
This system represents an extension of a test example
used in~\cite{lohmann-eid}, where 4~masses were included.

In the stochastic modeling, we replace the parameters by independent
uniformly distributed random variables, which vary 10\% around their
mean values.
The mean values are the constant parameters from above.
In the truncated orthogonal expansion~(\ref{truncated-pce}),
we include all multivariate Legendre polynomials up to degree three.
We obtain $m=1140$ basis functions.
The stochastic Galerkin system~(\ref{galerkin}) is calculated
exactly except for round-off errors.
Table~\ref{tab:msd} illustrates the properties of this example.
The system is asymptotically stable.
The mass matrix is symmetric and positive definite.
Yet the system is not dissipative.

\begin{table}
  \caption{Properties of stochastic Galerkin system in
    mass-spring-damper example.}
  \begin{center}
    \begin{tabular}{cc} \hline
      dimension & 11400 \\
      number of outputs & 1140 \\
      \# non-zeros in $\hat{E}$ & 13110 \\
      \# non-zeros in $\hat{A}$ & 50958 \\
      spectral abscissa of $(\hat{E},\hat{A})$ & $-0.0036$ \\ 
      spectral abscissa of $\hat{A} + \hat{A}^{\top}$ & 42.69 \\ \hline
    \end{tabular}
  \end{center}
\label{tab:msd}
\end{table}

We employ the one-sided Arnoldi method, see~\cite{antoulas},
to perform an MOR of the stochastic Galerkin system~(\ref{galerkin}).
The single real expansion point $s=0.7$
is used in this Krylov subspace technique.
The reduced systems are arranged for dimensions $r=1,\ldots,100$.
It follows that just 38 ROMs are stable.

Now we investigate the stabilization techniques from
Section~\ref{sec:stability}.
The following cases are discussed:
\begin{itemize}
   \item[i)] transformation of stochastic Galerkin system
    (Section~\ref{sec:trafo-gal}), 
   \item[ii)] transformation of original systems (Section~\ref{sec:trafo-org}),
    and 
    \item[iii)] transformation based on reference parameter
    (Section~\ref{sec:ref-parameter}). 
\end{itemize}

  In (i) and (iii), we reuse the projection matrix~$V$ determined
  for the original stochastic Galerkin system,
  because the Arnoldi method is invariant with respect to a
  basis transformation in the image space.
  In~(ii), we repeat the Arnoldi algorithm for the alternative
  stochastic Galerkin system, because a smaller reduction error
  is expected.
  
\begin{figure}
  \begin{center}
    untransformed \hspace{5cm} (i)

    \hspace{1mm}
    
  \includegraphics[width=6.5cm]{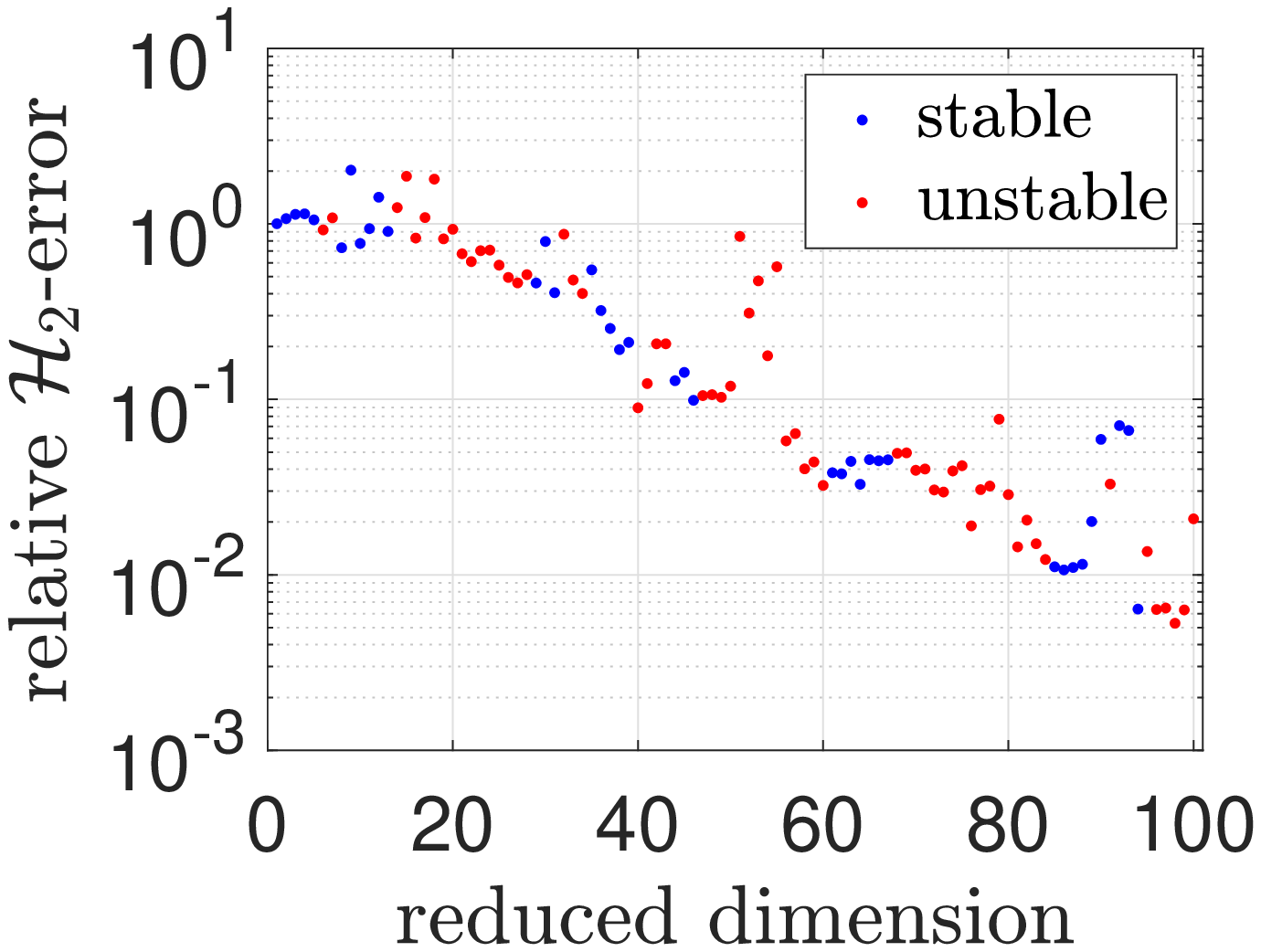}
  \hspace{5mm}
  \includegraphics[width=6.5cm]{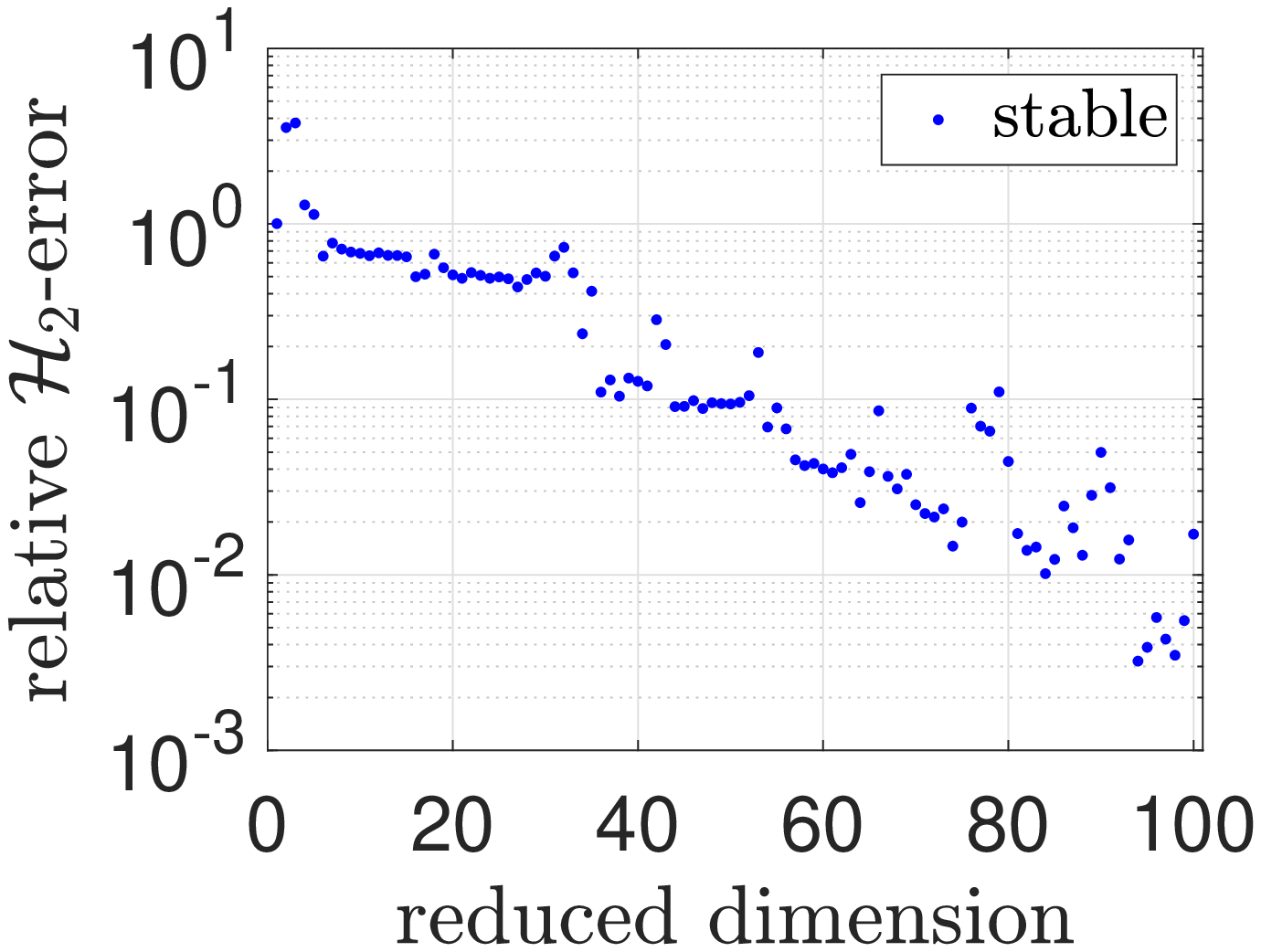}

  \hspace{5mm}

    (ii) \hspace{6.5cm} (iii)

    \hspace{1mm}
  
  \includegraphics[width=6.5cm]{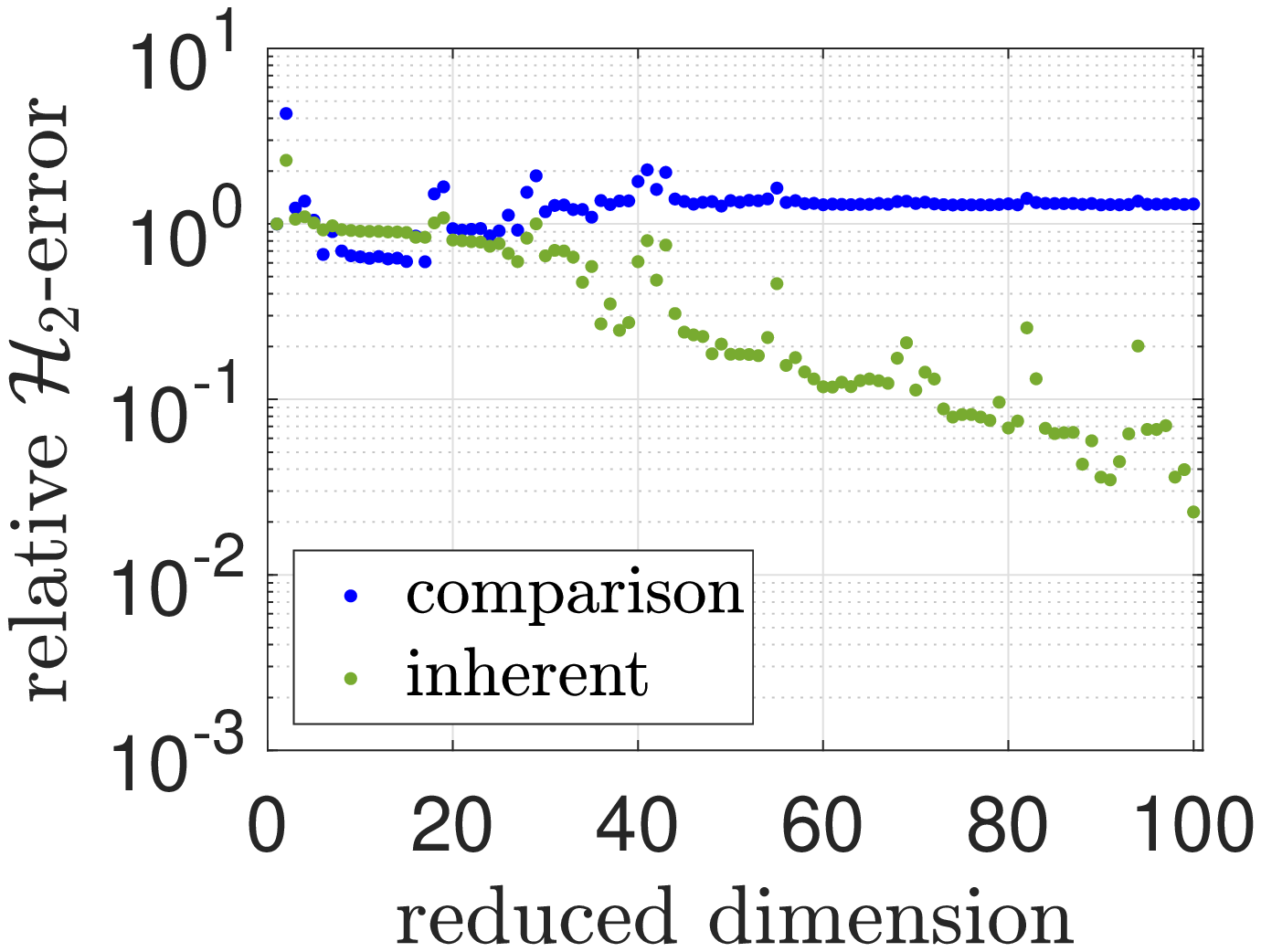}
  \hspace{5mm}
  \includegraphics[width=6.5cm]{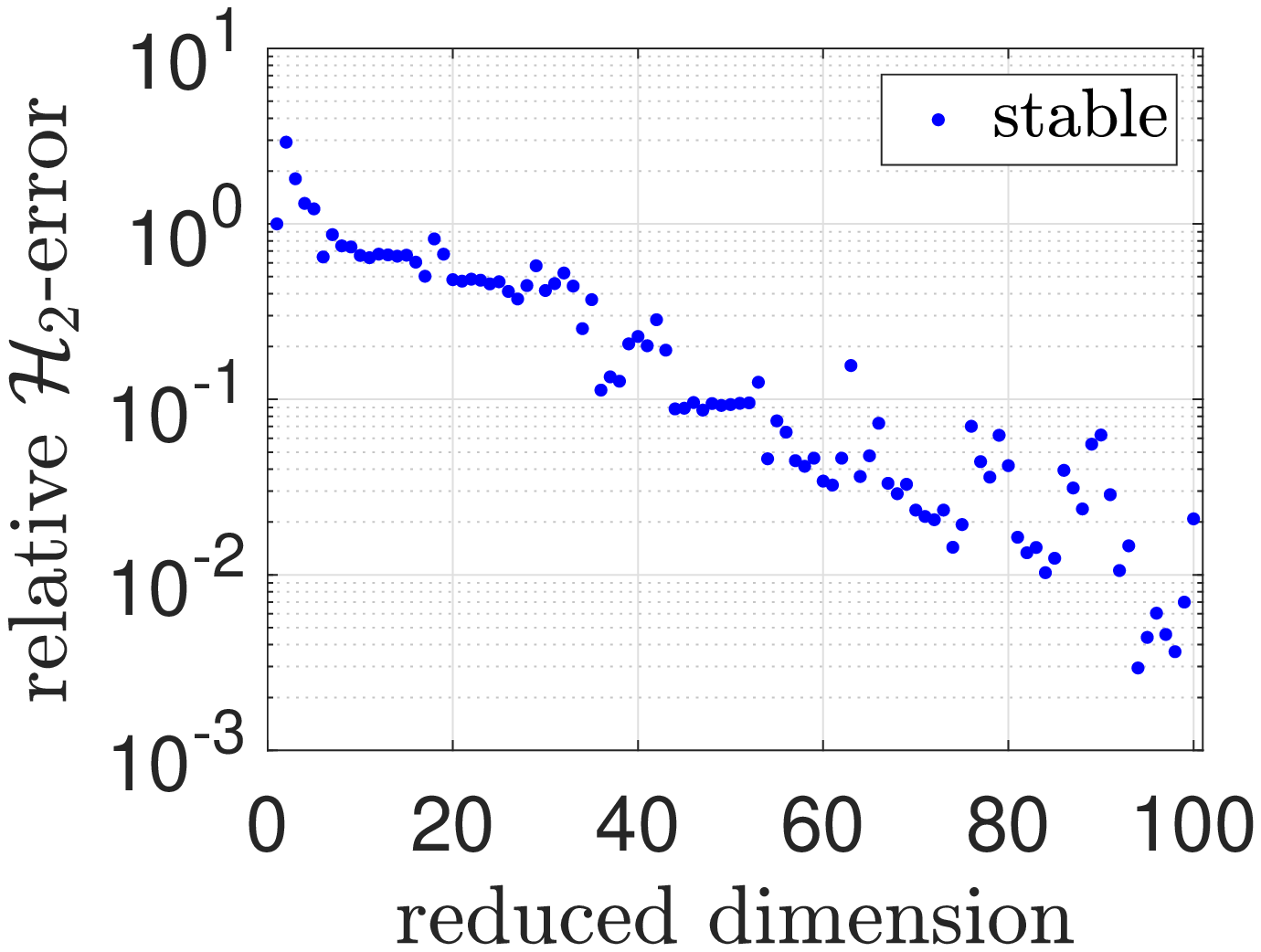}
  \end{center}
  \caption{Relative errors in $\htwo$-norm for MOR
    of stochastic Galerkin system and stabilization techniques
    in mass-spring-damper example.}
\label{fig:msd-errors}
\end{figure}

We choose the identity matrix ($F=I$) as degree of freedom in each
Lyapunov equation.
In all three approaches, the stability is achieved for each ROM.
Figure~\ref{fig:msd-errors} illustrates the relative errors
with respect to the $\htwo$-norm~(\ref{htwo-norm}), i.e.,
\begin{equation} \label{relative-error}
  \frac{\| H - \bar{H} \|_{\htwo}}{\| H \|_{\htwo}}
\end{equation}
with the transfer functions $H,\bar{H}$ of FOM and ROM, respectively.

\begin{figure}
  \begin{center}  
  \includegraphics[width=6.5cm]{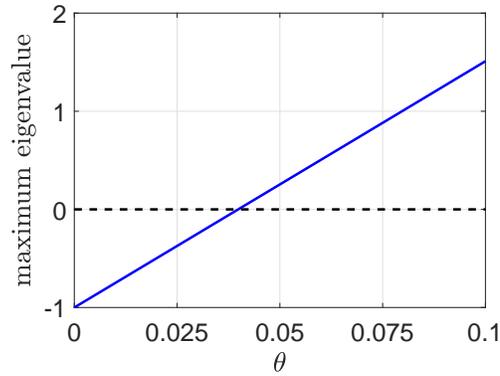}
  \end{center}
  \caption{Maximum eigenvalue of symmetric matrix
    $\hat{E}^\top \hat{M} \hat{A} + (\hat{E}^\top \hat{M} \hat{A})^\top$
    with $\hat{M}$ from~(\ref{large-M}) in technique (iii).}
\label{fig:msd-eigenvalues}
\end{figure}

\begin{table}
  \caption{Number of asymptotically stable systems
    out of 100 ROMs with dimensions $1,\ldots,100$
      and computation times for projection
      matrix~(\ref{projection-w})
    in stabilization using frequency domain integrals
    for mass-spring-damper example.}
  \begin{center}
  \begin{tabular}{lcccccccccc}
    \# nodes in quadrature & & 10 & & 20 & & 30 & & 40 & & untransf. \\ \hline
    \# stable ROMs & & 79 & & 91 & & 96 & & 100 & & 38 \\
    computation time (s) & & 19.0 & & 38.0 & & 56.8 & & 75.6 & & --
  \end{tabular}
  \end{center}
\label{tab:freq-integral}
\end{table}

\begin{table}
  \caption{Number of non-zero entries (nnz) in matrices of stochastic
    Galerkin systems
    for mass-spring-damper example.}
  \begin{center}
  \begin{tabular}{lcccc}
     & nnz in $\hat{E}$ & & nnz in $\hat{A}$ & \\ \hline
    untransformed Galerkin system &
    $13\,110$ & (0.01\%) & $50\,958$ & (0.04\%) \\
    Galerkin system in (ii) &
    $30\,000\,000$ & (23.1\%)& $29\,999\,562$ & (23.1\%)
  \end{tabular}
  \end{center}
\label{tab:sparsity}
\end{table}

In the technique~(i), we use the method from~\cite{pulch-arxiv},
where an integral is discretized by a quadrature rule in the frequency domain.
Therein, the solution of the high-dimensional Lyapunov equation is not
computed but the associated projection matrix~(\ref{projection-w}).
The computation work is characterized by an $LU$-decomposition of the
system matrix 
${\rm i}\omega_{j}\hat{E}-\hat{A}$ in each node~$\omega_{j}$
for the angular frequency~$\omega$. 
We employ the (univariate) Gauss-Legendre quadrature.
Table~\ref{tab:freq-integral} shows the number of stable ROMs for
different numbers of nodes.
We observe that 40~nodes are sufficient to stabilize all reduced systems.
Figure~\ref{fig:msd-errors} (i) depicts the relative
$\htwo$-errors~(\ref{relative-error}) in this case.
The error of the MOR exhibits the same magnitude as in the original
stochastic Galerkin system.

In the technique~(ii), we use a sparse grid quadrature of Smolyak-type
with level~3 based on the Clenshaw-Curtis rule.
The scheme exhibits 7209 nodes in the 17-dim\-ensional space
and negative weights arise.
Table~\ref{tab:sparsity} illustrates the loss of sparsity in this
approach.
Again all ROMs become asymptotically stable. 
The relative errors~(\ref{relative-error}) of the MOR are demonstrated in
Figure~\ref{fig:msd-errors} (ii).
The errors between the ROMs and the stochastic Galerkin projection
of the transformed systems, which we call the inherent errors, 
decay for increasing dimensions.
However, the errors between these ROMs and the stochastic Galerkin system
of the untransformed systems~(\ref{ode-parameter}),
which we call the comparison, are large and stagnate.
This property indicates that the quadrature error is too large,
even though a high number of nodes is used,
because the original Galerkin system can be considered
sufficiently accurate.

\newpage

In the technique~(iii), we define the mean value of the random variables
as reference parameter.
We calculate the matrices of the transformed Galerkin system to analyze
the definiteness.
Figure~\ref{fig:msd-eigenvalues} depicts the maximum eigenvalue of 
the symmetric part of $\hat{E}^\top \hat{M} \hat{A}$ with
$\hat{M}$ from~(\ref{large-M}) in dependence on the parameter~$\theta$
from~(\ref{random-variables}).
We observe that the negative definiteness is lost for
$\theta \ge \theta_0 \approx 0.04$.
Although our relevant case $\theta = 1$ is not dissipative,
the stability is preserved in all ROMs.
Moreover, the error of the MOR is not compromised.

\begin{table} 
  \caption{Computation times for MOR method producing reduced dimension
    $r=100$ and additional stabilization techniques
    in mass-spring-damper example.}
  \begin{center}
    \begin{tabular}{lc}
      & computation time (s) \\ \hline
      Arnoldi method for original system & 2.4 \\
      Arnoldi method for  system in~(ii) & 65.4 \\
      technique (i) & 75.6 \\
      technique (ii) & 2591.9 \\
      technique (iii) & 0.1 
    \end{tabular}
  \end{center}
\label{tab:computation-times}
\end{table}

We also measured computation times illustrated by
Table~\ref{tab:computation-times}.
The ROMs of dimension $r=100$ are determined using the Arnoldi method
and the three stability-preserving techniques.
The computation time of a stabilization technique consists of all steps
to calculate the reduced matrices~(\ref{matrices-reduced})
except for the determination of the projection matrix~$V$,
which is done by the Arnoldi method.
In the technique~(i), $k=40$ nodes are used in the Gaussian quadrature.
On the one hand, the technique~(ii) causes a relatively large computational
effort due to the high number of nodes in the sparse grid quadrature.
Also the Arnoldi algorithm is more expensive due to the loss of sparsity
in the matrices. 
On the other hand, the technique~(iii) requires a negligible computation
work, which shows that it is worth trying this approach.

\subsection{Band-pass filter}
\begin{figure}[h!]
  \begin{center}  
  \includegraphics[width=13.5cm]{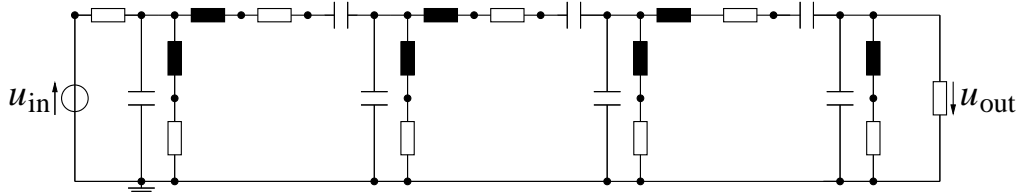}
  \end{center}
  \caption{Electric circuit of band-pass filter
    (white boxes: resistors, black boxes: inductors,
    parallel lines: capacitors).}
\label{fig:filter-circuit}
\end{figure}

\begin{figure}
  \begin{center}  
  \includegraphics[width=6.5cm]{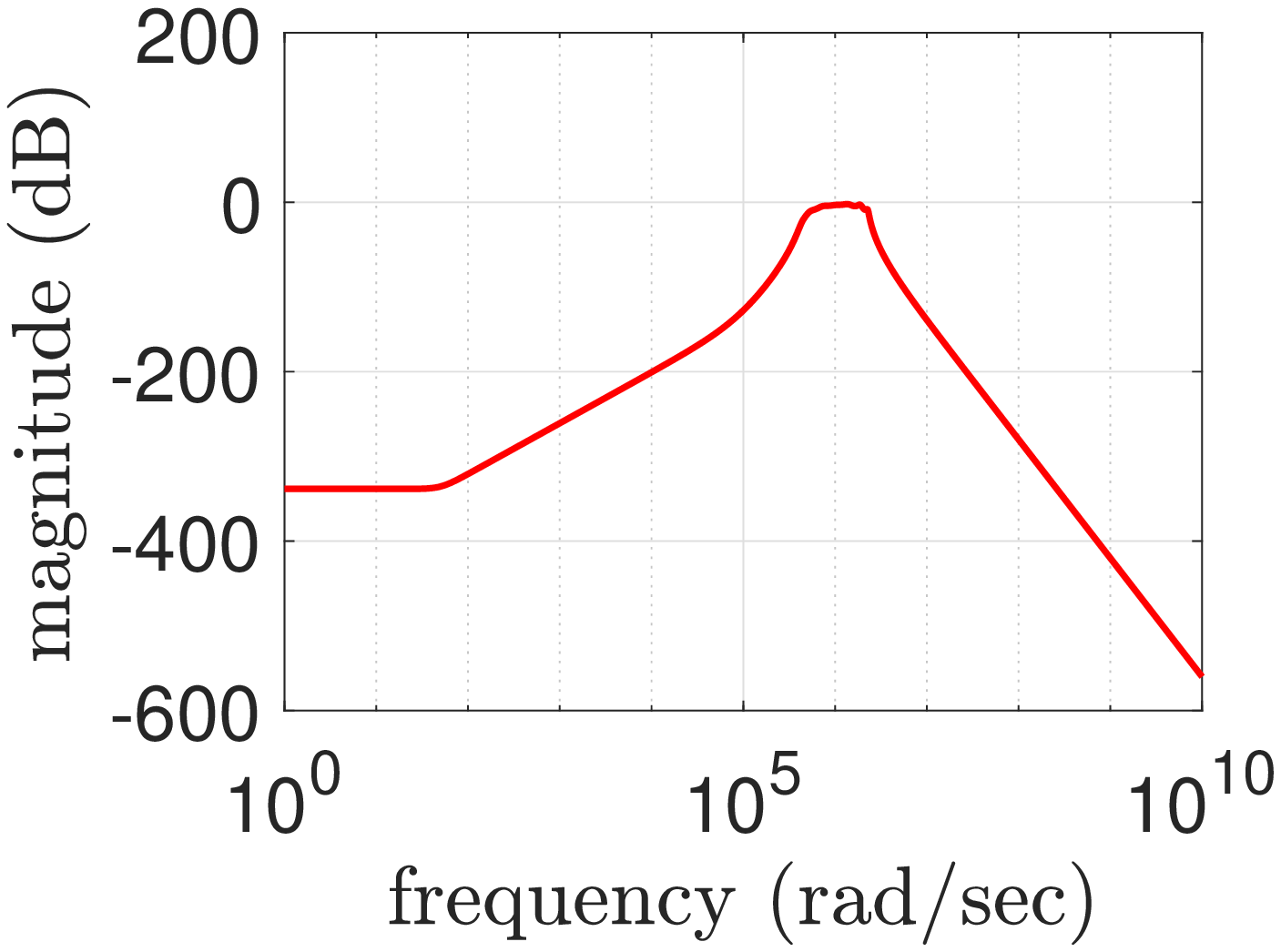}
  \hspace{5mm}
  \includegraphics[width=6.5cm]{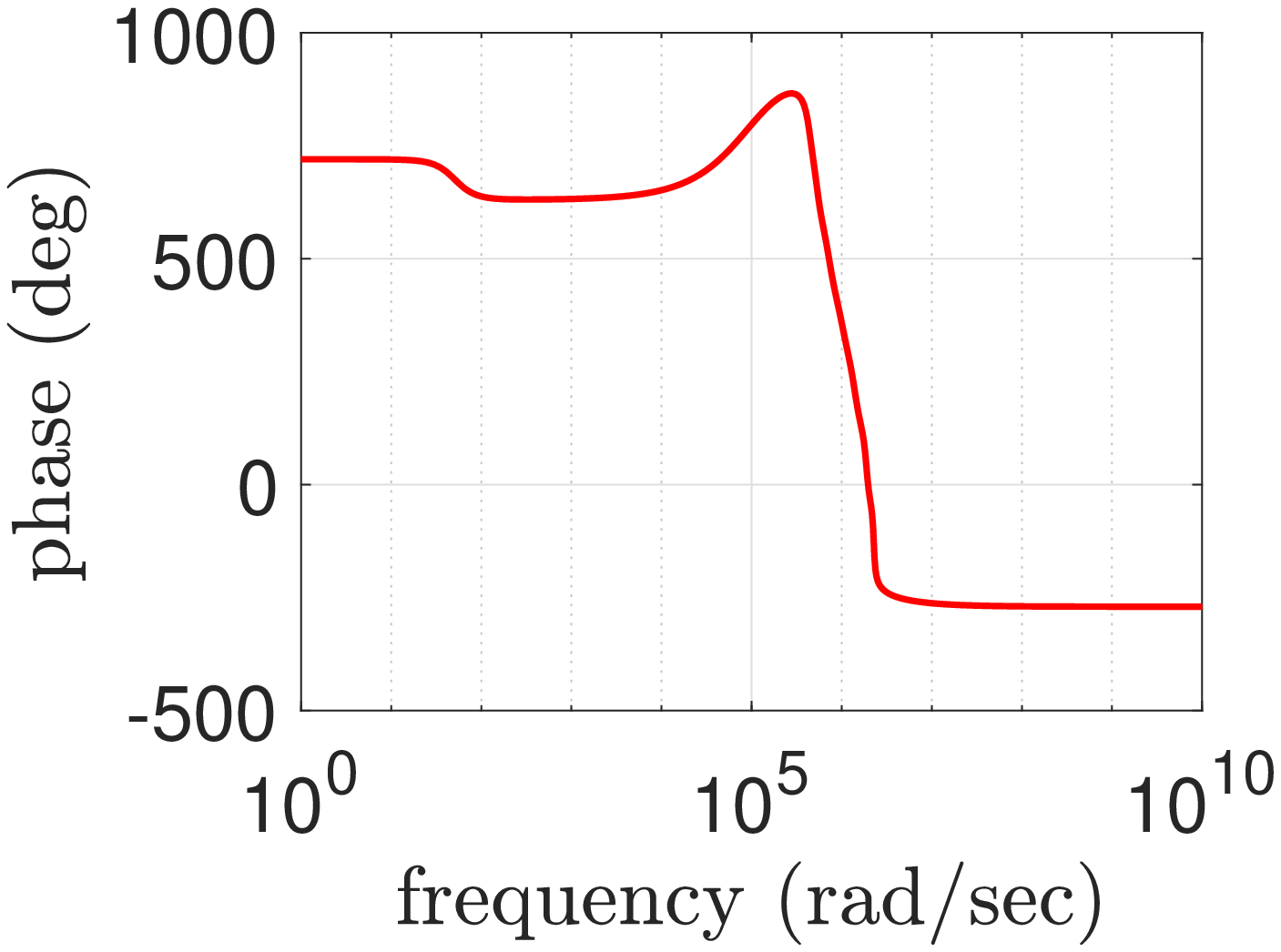}
  \end{center}
  \caption{Bode plot of band-pass filter for
  a constant choice of the parameters.}
\label{fig:bpf-bode}
\end{figure}

We examine the electric circuit of a band-pass filter
shown in Figure~\ref{fig:filter-circuit}.
A single input voltage is supplied and a single output voltage
drops at a load conductance.
Modified nodal analysis~\cite{ho} yields a linear system of DAEs of
dimension $n=23$.
This DAE system exhibits the index one.
The physical parameters are 7~capacitances, 7~inductances
and 9~conductances ($q=23$).
Figure~\ref{fig:bpf-bode} depicts the Bode plot of the system for
a constant selection of the parameters.
Furthermore, the system matrices are affine-linear functions
of the parameters. 

\begin{table}[t]
  \caption{Properties of regularized stochastic Galerkin system in
    band-pass filter example.}
  \begin{center}
    \begin{tabular}{cc} \hline
      dimension & 6900 \\
      number of outputs & 300 \\
      \# non-zeros in $\hat{E}_{\rm reg}$ & 21912 \\
      \# non-zeros in $\hat{A}_{\rm reg}$ & 21912 \\
      spectral abscissa of $(\hat{E}_{\rm reg},\hat{A}_{\rm reg})$
      & $-4.94 \cdot 10^{-4}$ \\ 
      spectral abscissa of
      $(\hat{E}_{\rm reg}^\top \hat{A}_{\rm reg}) +
      (\hat{E}_{\rm reg}^\top \hat{A}_{\rm reg})^\top$
      & $1.28 \cdot 10^{-5}$ \\ \hline
    \end{tabular}
  \end{center}
\label{tab:bpf}
\end{table}

In the stochastic modeling, we introduce independent uniformly distributed
random variables varying 20\% around their mean values.
The truncated series~(\ref{truncated-pce}) include all basis polynomials
up to degree two, i.e., $m=300$.
The stochastic Galerkin system~(\ref{galerkin}) consists of linear DAEs,
which have a finite $\htwo$-norm~(\ref{htwo-norm}) for the defined outputs.

In the regularization~(\ref{regularization}),
we choose the parameters $\alpha = 10^{-10}$ and $\beta = 10^{-5}$.
Table~\ref{tab:bpf} shows the properties of the regularized system.
The system is asymptotically stable. 
Both $\hat{E}$ and $\hat{E}_{\rm reg}$ are not symmetric.
Thus the spectral abscissa, see Definition~\ref{def:abscissa},
of the symmetric part of $\hat{A}_{\rm reg}$ is irrelevant.
Multiplication by $\hat{E}_{\rm reg}^\top$ from the left yields a system
with symmetric positive definite mass matrix,
which is still not dissipative.

The one-sided Arnoldi method yields the projection matrices of the MOR
for dimensions $r=1,\ldots,100$.
Therein, we choose the real number $s=10^6$ as expansion point. 
A loss of stability occurs in the reduction of both the
DAE system and the regularized system,
as shown in Table~\ref{tab:bpf-stable}.
We apply the stabilization techniques of
Section~\ref{sec:trafo-gal} and Section~\ref{sec:ref-parameter}
to the regularized system.
The Lyapunov equations are solved using the identity matrix as input matrix.
We solve the Lyapunov equations by a direct method of linear algebra.
Thus a critical behavior for small regularization parameters,
as mentioned in Section~\ref{sec:dae-gal}, is avoided.
The matrices of a transformed Galerkin system are never
computed explicitly.
All ROMs become asymptotically stable in both approaches.

\begin{figure}
  \begin{center}
    (i) \hspace{6.5cm} (ii)

    \hspace{1mm}
    
  \includegraphics[width=6.5cm]{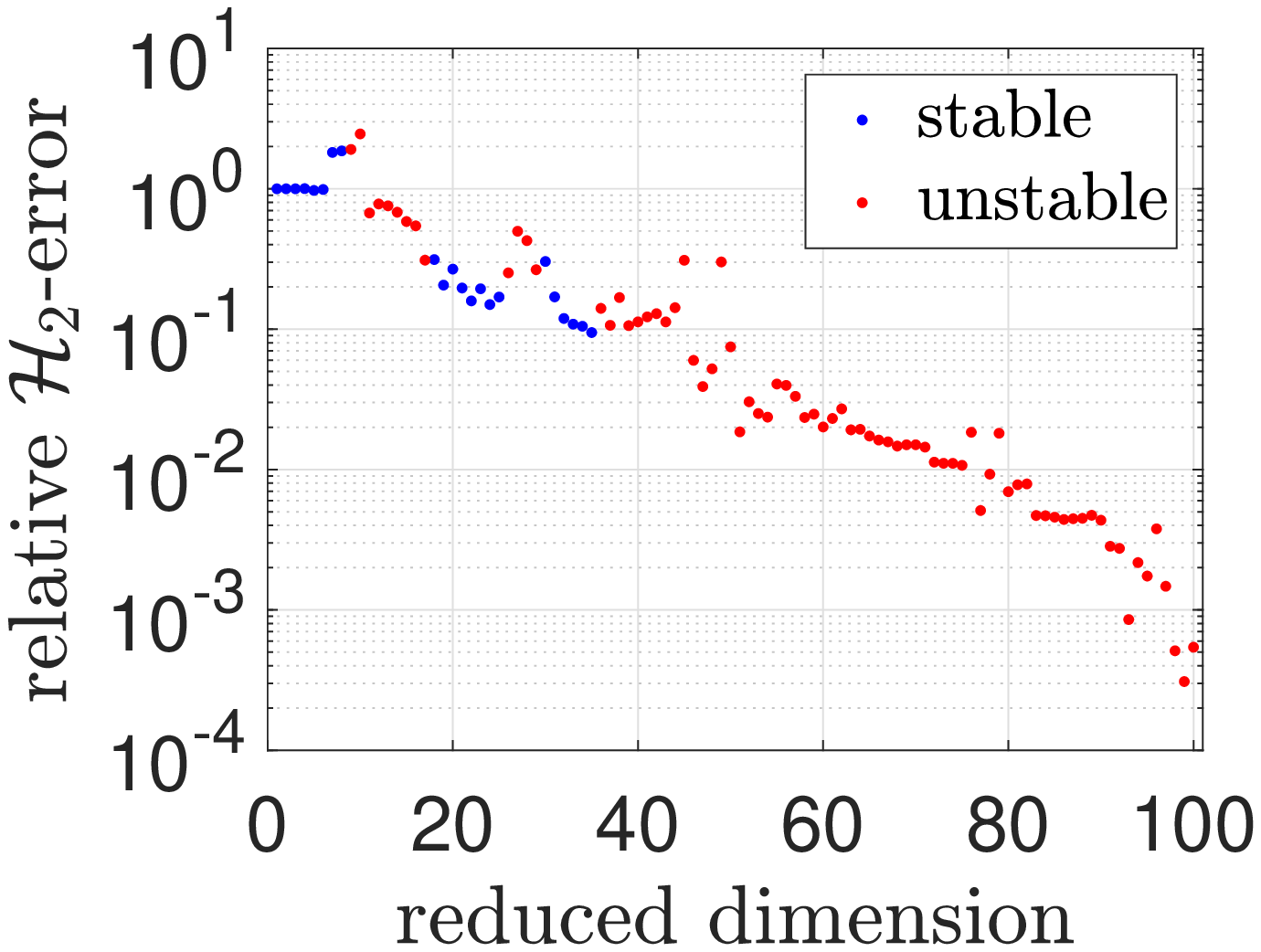}
  \hspace{5mm}
  \includegraphics[width=6.5cm]{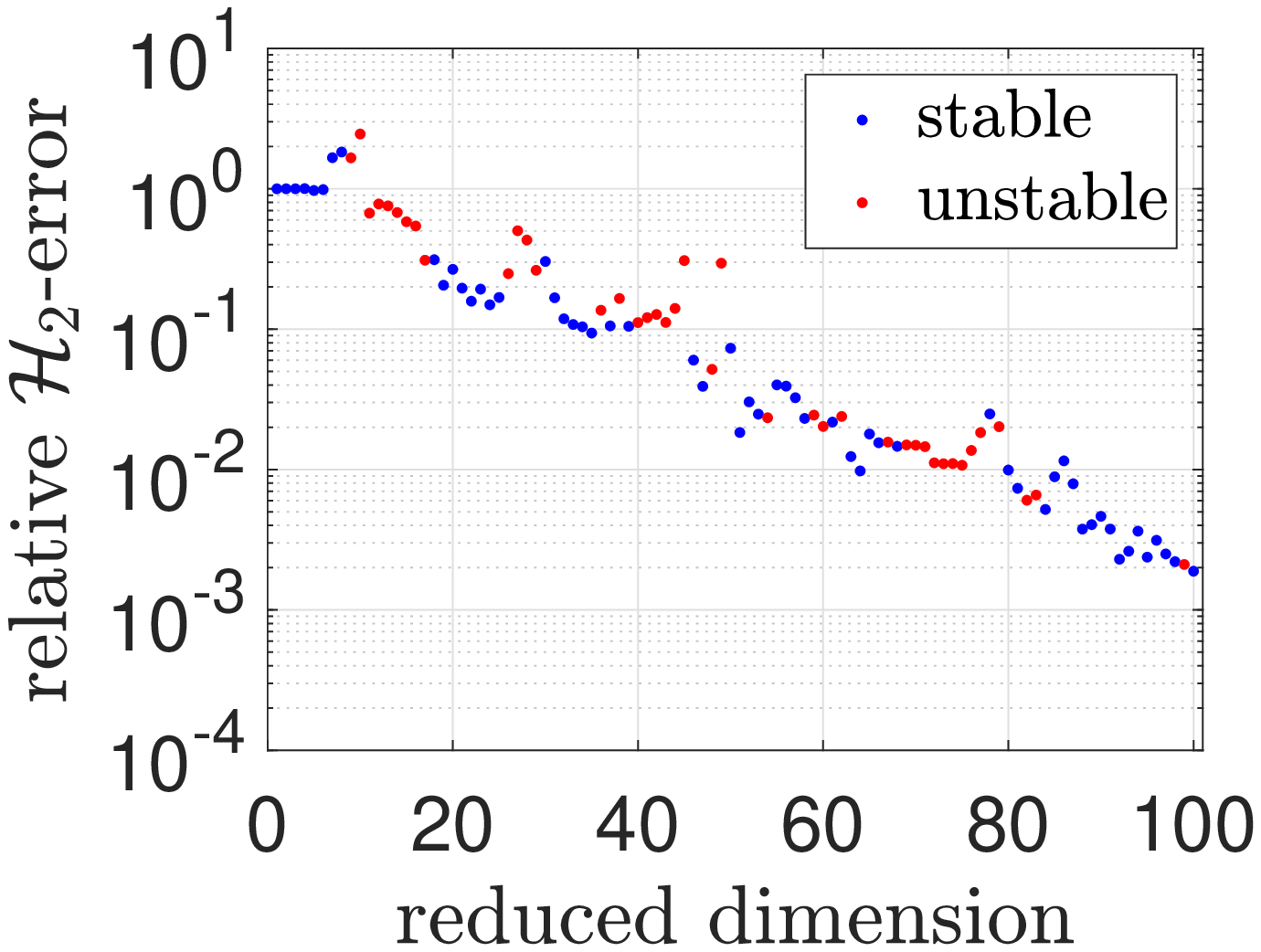}

  \hspace{5mm}

    (iii) \hspace{6.5cm} (iv)

    \hspace{1mm}
  
  \includegraphics[width=6.5cm]{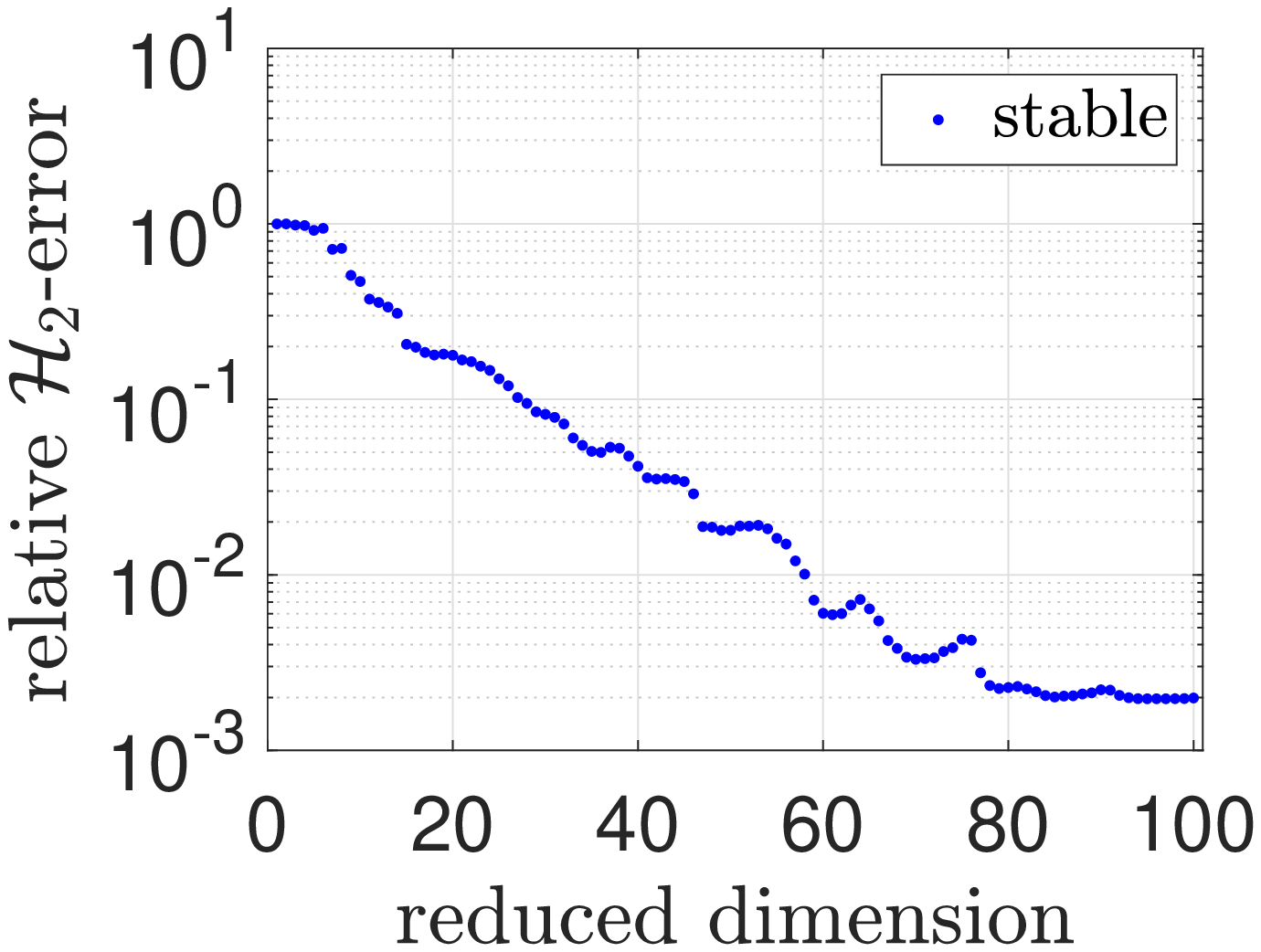}
  \hspace{5mm}
  \includegraphics[width=6.5cm]{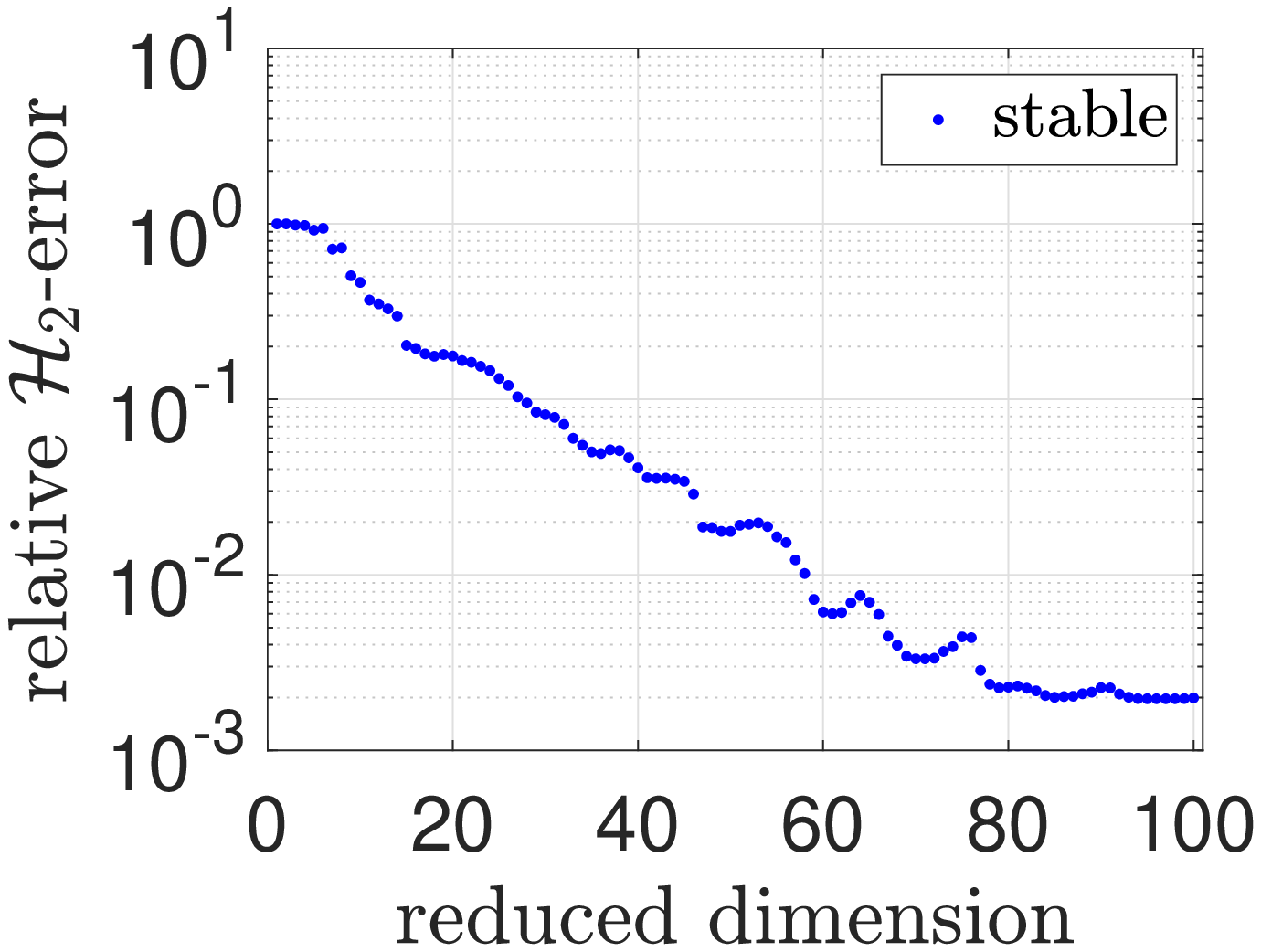}
  \end{center}
  \caption{Relative errors in $\htwo$-norm for MOR
    of different FOMs and stabilization techniques,
    see Table~\ref{tab:bpf-stable},
    in band-pass filter example.}
\label{fig:bpf-errors}
\end{figure}

Finally, we compute approximations of the $\htwo$-norms for
the difference between the DAE system and the reduced systems.
Figure~\ref{fig:bpf-errors} illustrates the
relative errors~(\ref{relative-error}) of the reductions.
We recognize that the errors are nearly identical in the two
stabilization techniques.
The errors stagnate for reduced dimensions $r>80$ in (ii)-(iv),
because the total error is dominated by the error of the
regularization in this part.
Most important, the stabilization approaches do not compromise
the error of the MOR.

\begin{table}
  \caption{Number of asymptotically stable systems
    out of 100 ROMs with dimensions $1,\ldots,100$ 
      in MOR of stochastic Galerkin system for band-pass filter example.}
  \begin{center}
    \begin{tabular}{rlc}
      & FOM system & \# stable ROMs \\ \hline
      i) & DAE & 22 \\
      ii) & regularized system (ODE) & 59 \\
      iii) & transformed as in Section~\ref{sec:trafo-gal} (ODE) & 100 \\
      iv) & transformed as in Section~\ref{sec:ref-parameter} (ODE) & 100
    \end{tabular}
  \end{center}
\label{tab:bpf-stable}
\end{table}

\clearpage


\section{Conclusions}
We examined stability-preserving model order reduction of
linear stochastic Ga\-lerkin systems using transformations to
dissipative forms.
Three approaches were analyzed.
The transformation of the conventional Galerkin system represents
an adequate method.
The Galerkin projection of the transformed original systems
features severe drawbacks in the case of large numbers of random parameters.
The approach using a cheap transformation matrix based on a
reference parameter is promising.
Although the dissipativity property cannot be guaranteed in the
transformed system, the numerical results of test examples demonstrate
that preservation of stability is achieved.
Moreover, the error of the model order reduction does not increase
in this stabilization.



\begin{thebibliography}{00}

\bibitem{antoulas}
Antoulas,~A.:
Approximation of Large-Scale Dynamical Systems.
SIAM Publications, 2005.

\bibitem{benner}
Benner, P.; Hinze, M.; ter Maten, E.J.W. (eds.):
Model Reduction for Circuit Simulation.
Lect. Notes in Electr. Engng. Vol.~74, Springer, 2011.

\bibitem{benner-gugercin}
  Benner,~P.; Gugercin,~S.; Willcox,~K.:
  A survey of projection-based model order reduction methods for
  parametric dynamical systems.
  SIAM Review~57 (2015) 483--531.

\bibitem{benner-mehrmann}
  Benner,~P.; Mehrmann,~V.; Sorensen,~D.C. (eds.):
  Dimension Reduction of Large-Scale Systems.
  Lect. Notes in Comput. Sci. Eng. Vol.~45,
  Springer, 2005.

\bibitem{benner-schneider}
  Benner,~P.; Schneider,~A.:
  Balanced truncation model order reduction for LTI systems
  with many inputs or outputs.
  In: Edelmayer, A. (eds.),
  Proceedings 19th International Symposium on Mathematical Theory
  of Networks and Systems, 2010, pp.~1971--1974.

\bibitem{benner-stykel}
  Benner,~P.; Stykel,~T.: 
  Model order reduction of differential-algebraic equations: a survey.
  In: Ilchmann,~A.; Reis,~T. (eds.),
  {Surveys in Differential-Algebraic Equations IV}, 
  Differential-Algebraic Equations Forum, Springer, 2017, pp. 107--160.

\bibitem{castane-selga}
  Casta{\~n}{\'e}~Selga,~R.; Lohmann,~B.; Eid,~R.:
  Stability preservation in projec\-tion-based model order reduction
  of large scale systems.
  Eur. J. Control 18 (2012) 122--132.

\bibitem{ernst-etal}
  Ernst,~O.G; Mugler,~A.; Starkloff,~H.J., Ullmann,~E.:
  On the convergence of generalized polynomial chaos expansions.
  ESAIM: M2AN 46 (2012) 317--339.

\bibitem{freitas}
  Freitas,~F.D.; Pulch,~R.; Rommes,~J.:
  Fast and accurate model reduction for spectral methods
  in uncertainty quantification.
  Int. J. Uncertain. Quantif. 6 (2016) 271--286.

\bibitem{freund}
  Freund,~R.:
  Model reduction methods based on Krylov subspaces.
  Acta Numerica 12 (2003) 267--319.

\bibitem{hammarling}
  Hammarling,~S.J.:
  Numerical solution of stable non-negative definite Lyapunov equation.
  IMA J. Numer. Anal. 2 (1982) 303--323.

\bibitem{ho}
  Ho,~C.W.; Ruehli,~A.; Brennan,~P.:
  The modified nodal approach to network analysis.
  IEEE Trans. Circ. Syst.~22 (1975) 504--509.
  
\bibitem{kramer-singler}
  Kramer,~B.; Singler,~J.R.:
  A POD projection method for large-scale algebraic Riccati equations.
  Numer. Algebra Contr. Optim. 6 (2016) 413--435.
  
\bibitem{lohmann-eid}
  Lohmann,~B.; Eid,~R.:
  Efficient order reduction of parametric and nonlinear models
  by superposition of locally reduced models.
  In: Roppenecker,~G.; Lohmann,~B. (eds.),
  Methoden und Anwendungen der Regelungstechnik, Shaker, 2009.

\bibitem{lu-wachspress}
  Lu,~A.; Wachspress,~E.L.:
  Solution of Lyapunov equations by alternating direction implicit iteration.
  Computers Math. Applic. 21 (1991) 43--58.

\bibitem{manfredi}
  Manfredi,~P.; Vande~Ginste,~D.; De~Zutter,~D.; Canavero, F.G.:
  On the passivity of polynomial chaos-based augmented models for
  stochastic circuits.
  IEEE Trans. Circuits Syst. I-Regul. Pap. 60 (2013) 2998--3007.
  
\bibitem{maten-nanocops}
  ter Maten, E.J.W. et al.:
  Nanoelectronic coupled problems solutions -- nanoCOPS:
  modelling, multirate, model order reduction, uncertainty quantification,
  fast fault simulation.
  J. Math. Ind.~7 (2016).
  
\bibitem{matlab2018}
  MATLAB, version 9.4.0.813654 (R2018a), The Mathworks Inc.,
  Natick, Massachusetts, 2018.

\bibitem{mi-etal}
  Mi,~N.; Tan,~S.X.D.; Liu,~P.; Cui,~J.; Cai,~Y.; Hong,~X.:
  Stochastic extended Krylov subspace method for variational
  analysis of on-chip power grid networks.
  Proc. ICCAD 2007, pp.~48--53.

\bibitem{mueller}
  M\"uller,~P.C.:
  Modified Lyapunov equations for LTI descriptor systems.
  J.~Braz. Soc. Mech. Sci. \& Eng. 28 (2006) 448--452.

\bibitem{panzeretal}
  Panzer,~H.; Wolf,~T.; Lohmann,~B.:
  A strictly dissipative state space representation of second order systems.
  Automatisierungstechnik 60 (2012) 392--396.

\bibitem{penzl1998}
  Penzl, T.:
  Numerical solution of generalized Lyapunov equations.
  Adv. Comput. Math. 8 (1998) 33--48.

\bibitem{penzl-sisc}
  Penzl,~T.:
  A cyclic low-rank Smith method for large sparse Lyapunov equations.
  SIAM J. Sci. Comput. 21 (2000) 1401--1418.

\bibitem{prajna}
  Prajna,~S.:
  POD model reduction with stability guarantee.
  In: Proceedings of 42nd IEEE Conference on Decision and Control,
  Maui, Hawaii, USA, 2003, pp.~5254--5258.

\bibitem{pulch14}
Pulch,~R.:
Stochastic collocation and stochastic Galerkin methods for 
linear differential algebraic equations.
J. Comput. Appl. Math. 262 (2014) 281--291.

\bibitem{pulch-scee2014}
Pulch,~R.:
Model order reduction for stochastic expansions of electric circuits.
In: Bartel,~A.; Clemens,~M.; G\"unther,~M.; ter Maten,~E.J.W. (eds.), 
Scientific Computing in Electrical Engineering SCEE 2014,
Mathematics in Industry Vol. 23, Springer, 2016, pp. 223--232.

\bibitem{pulch17}
Pulch,~R.:
Model order reduction and low-dimensional representations for
random linear dynamical systems.
Math. Comput. Simulat. 144 (2017) 1--20.

\bibitem{pulch-arxiv}
  Pulch,~R.:
  Frequency domain integrals for stability preservation in
  Galerkin-type pro\-jec\-tion-based model order reduction.
  {\tt arxiv:1808.04119} (2018).

\bibitem{pulch-naco}
  Pulch,~R.:
  Stability preservation in Galerkin-type projection-based
  model order reduction.
  Numer. Algebra Contr. Optim. 9 (2019) 23--44.

\bibitem{pulch19}
  Pulch,~R.:
  Model order reduction for random nonlinear dynamical systems and
  low-dimensional representations for their quantities of interest.
  Math. Comput. Simulat. 166 (2019) 76--92.

\bibitem{pulch-augustin}
  Pulch,~R.; Augustin,~F.:
  Stability preservation in stochastic Galerkin projections
  of dynamical systems.
  SIAM/ASA J. Uncertainty Quantification 7 (2019) 634--651.

\bibitem{pulch-maten}
Pulch,~R.; ter Maten,~E.J.W.: 
Stochastic Galerkin methods and model order reduction for
linear dynamical systems. 
Int. J. Uncertain. Quantif. 5 (2015) 255--273.

\bibitem{schilders}
Schilders,~W.H.A.; van~der~Vorst,~M.A.; Rommes,~J. (eds.): 
Model Order Reduction: Theory, Research Aspects and Applications. 
{Mathematics in Industry} Vol.~13, Springer, 2008.

\bibitem{seydel}
  Seydel,~R.:
  Practical Bifurcation and Stability Analysis. (3rd ed.) Springer, 2010.
  
\bibitem{son-stykel}
  Son,~N.T.; Stykel,~T.:
  Solving parameter-dependent Lyapunov equations using the reduced basis
  method with application to parametric model order reduction.
  SIAM J. Matrix Anal. Appl. 38 (2017) 478--504.

\bibitem{sonday}
  Sonday,~B.; Berry,~R.; Debusschere,~B.; Najm,~H.:
  Eigenvalues of the Jacobian of a Galerkin-projected uncertain ODE system.
  SIAM J. Sci. Comput.~33 (2011) 1212--1233.

\bibitem{sullivan}
  Sullivan, T.J.:
  Introduction to Uncertainty Quantification.
  Springer, 2015.
  
\bibitem{wolf}
  Wolf,~T.; Panzer,~H.; Lohmann,~B.:
  Model order reduction by approximate balanced truncation:
  a unifying framework.
  Automatisierungstechnik~61 (2013) 545--556.
  
\bibitem{xiu-book}
  Xiu,~D.:
  Numerical methods for stochastic computations: a spectral method approach.
  Princeton University Press, 2010.

\bibitem{zou-etal}
  Zou,~Y.; Cai,~Y.; Zhou,~Q.; Hong,~X.; Tan,~S.X.-D.; Kang,~L.:
  Practical implementation of the stochastic parameterized
  model order reduction via Hermite polynomial chaos.
  Proc. ASP-DAC 2007, pp. 367-–372.
\end{thebibliography}
\end{document}